\documentclass[12pt]{article}

\newtheorem{theo}{Theorem}
\newtheorem{prop}{Proposition}

\newtheorem{lemma}{Lemma}
\newtheorem{remark}{Remark}

\newcommand{\Ref}[1]{(\ref{#1})}

\newenvironment{proof}{\textbf{Proof. }}
{$\bigtriangleup$}

\newenvironment{proof_theo1}{\noindent\textbf{Proof of theorem 1. }}
{$\bigtriangleup$}

\newenvironment{proof_theo2}{\noindent\textbf{Proof of theorem 2. }}
{$\bigtriangleup$}

\newenvironment{proof_theo3}{\noindent\textbf{Proof of theorem 3. }}
{$\bigtriangleup$}

\newenvironment{proof_theo4}{\noindent\textbf{Proof of theorem 4. }}
{$\bigtriangleup$}

\newenvironment{eq}{\begin{equation}}{\end{equation}}

\newcommand{\si}{\sigma}

\newcommand{\al}{\alpha}
\newcommand{\be}{\beta}
\newcommand{\ga}{\gamma}
\newcommand{\la}{\lambda}
\newcommand{\de}{\delta}
\newcommand{\De}{\Delta}

\newcommand{\LA}{\langle}
\newcommand{\RA}{\rangle}
\newcommand{\ov}[1]{\overline{#1}}

\newcommand{\tr}{{\rm tr}}
\newcommand{\sys}{{\mathcal{S}} }

\newcommand{\mdeg}{\mathop{\rm mdeg}}

\newcommand{\mod}[1]{|\,#1|}
\newcommand{\lin}{\mathop{\rm lin}}

\newcommand{\MSG}{G}

\newcommand{\un}[1]{{\underline{#1}} }

\newcommand{\OK}[1]{N_{3,d}(#1)}

\begin{document}
 \title{Relatively Free Algebras with the Identity $x^3=0$}

  \author{A. A. Lopatin
  \\ Chair of Algebra,
 \\ Department of Mathematics, \\
 Omsk State University,\\ 55A Prospect Mira, Omsk 644077, Russia,
 \\e-mail: lopatin@math.omsu.omskreg.ru }
 \date{}
\maketitle

\begin{abstract}
A basis for a relatively free associative algebra with the
identity $x^3=0$ over a field of an arbitrary characteristic is
found. As an application, a minimal generating system  for the
$3\times3$ matrix invariant algebra is determined.
\end{abstract}
\section{Introduction}
Let $K$ be an infinite field of an arbitrary characteristic $p$
($p=0,2,3,\ldots$). Let $K\LA x_1,\ldots,x_d\RA^{\#}$ ($K\LA
x_1,\ldots,x_d\RA$, respectively) be the free associative
$K$-algebra without unity (with unity, respectively) which is
freely generated by $x_1,\ldots,x_d$. Let ${\rm
id}\{f_1,\ldots,f_s\}$ be the ideal generated by $f_1,\ldots,f_s$.
Denote by $N_{n,d}=K\LA x_1,\ldots,x_d\RA^{\#}/$ ${\rm
id}\{x^n|x\in K\LA x_1,\ldots,x_d\RA^{\#}\}$ a relatively free
finitely generated non-unitary $K$-algebra with the identity
$x^n=0$, where $d\geq1$. Let ${\mathcal N}=\{0,1,2,\ldots\}$. The
algebra $N_{n,d}$ possesses natural ${\mathcal N}$- and ${\mathcal
N}^d$-gradings by degrees and multidegrees  respectively.

The {\it nilpotency degree} of a non-unitary algebra $A$ is the
least $C>0$ for which $a_1\cdots a_C=0$ for all $a_1,\ldots,a_C\in
A$. Denote by $C(n,d,K)$ the nilpotency degree of $N_{n,d}$. In
the case of characteristic zero $n(n+1)/2\leq C(n,d,K)\leq n^2$
(see~\cite{Kuzmin},~\cite{Raz_book}), and there is a conjecture
that $C(n,d,K)=n(n+1)/2$. This conjecture has been proven for
$n\leq4$ (see~\cite{Vau}). If $p=0$ or $p>n$, then $C(n,d,K)<2^n$
by~\cite{Higman}. For a positive characteristic some upper bounds
on $C(n,d,K)$ are given in~\cite{Klein}: $C(n,d,K)<(1/6)n^6d^n$
and $C(n,d,K)<1/(m-1)!\,n^{n^3}d^m$, where $m=[n/2]$.
In~\cite{Lopatin} $C(3,d,K)$ was established for an arbitrary
$d,p$, except for the case of $p=3$, $d$ is odd, where the
deviation in the estimation of $C(3,d,K)$ is equal to $1$. In this
article, a basis for $N_{3,d}$ is found (see
Proposition~\ref{prop_p0} and Theorems~\ref{theo2},~\ref{theo3}),
and, in particular, $C(3,d,K)$ is established for any $d,p$.
Namely, when $d>1$, we have:
$$
\begin{array}{l}
\mbox{If } p=0 \mbox{ or } p>3, \mbox{ then }C(3,d,K)=6.\\
\mbox{If }$p=2$, \mbox{ then }C(3,d,K)=\left\{
\begin{array}{lcl}
d+3&,& d\geq3\\
6&,& d=2.\\
\end{array}
\right. \\
\mbox{If }p=3, \mbox{ then }C(3,d,K)=3d+1.\\
\end{array}
$$

As an application, a minimal homogeneous generating system  of the
$3\times3$ matrix invariant algebra is determined (see
Theorem~\ref{theo4}).

For $p=2,3$, a basis for the multilinear homogeneous component of
$N_{3,d}$ for 'small' $d$ was found by means of a computer
programme. Then, the case of an arbitrary $d$ was reduced to the
case of 'small' $d$ using the composition method. All programmes
were written by means of Borland C++ Builder (version 6.0) and are
available upon request from the author. The notion of the
composition method was taken from~\cite{Bokut}.

\section{Preliminaries}

Further, we assume that $n=3$, unless it is stated otherwise. Let
${\mathcal Z}$ be the ring of integers, and let ${\mathcal Q}$ be
the field of rational fractions. Denote by $F_d$ the free
semigroup, generated by letters $\{x_1,x_2,\ldots,x_d\}$. By
$F_d^{\#}$ we mean $F_d$ without unity. For short, we will write
$K\LA F_d\RA$ instead of $K\LA x_1,\ldots, x_d\RA$. The degree of
a ${\mathcal N}^d$-homogeneous element $u\in F_d$ we denote by
$\deg(u)$, its multidegree we denote by $\mdeg(u)$, and the degree
of $u$ in letter $x_j$ we denote by $\deg_{x_j}(u)$. Elements of
$F_d$ are called words. By words from $N_{3,d}$ we mean images of
words from $F_d$ in $N_{3,d}$ under the natural homomorphism. We
assume that all words are non-empty, that is they are not equal to
unity of $F_d$, unless it is stated otherwise. Notation
$w=x_{i_1}\cdots \widetilde{x}_{i_s}\cdots x_{i_t}$ stands for the
word $w$, which can be get from the word $x_{i_1}\cdots x_{i_t}$
by eliminating the letter $x_{i_s}$. For a set of words $M$ and a
word $v$ denote by $vM$ the set $\{vu|\, u\in M\}$. If the set $M$
is empty, then we assume $vM=\emptyset$. By $x_i\cdots x_j$
($i,j\in \mathcal Z$) we mean the word $x_ix_{i+1}x_{i+2}\cdots
x_j$ if $1\leq i\leq j$, and the empty word otherwise.

For some ${\mathcal N}^d$-graded algebra $A$ and multidegree
$\De=(\de_1,\ldots,\de_d)$ denote by $A(\De)$ or
$A(\de_1,\ldots,\de_d)$ the homogeneous component of $A$ of
multidegree $\De$. For short,  multidegree
$(3,\ldots,3,2,\ldots,2,1\ldots,1)$ will be denoted by
$3^{r}2^{s}1^{t}$ for the appropriate $r,s,t$. For
$\De=(\de_1,\ldots,\de_t)$ let $\mod{\De}=\sum_i \de_i$. By
$\lin\{v_1,\ldots,v_t\}$ we mean the linear span of the elements
$v_1,\ldots,v_t$ of some vector space over $K$. We denote some
elements of $K\LA F_d\RA^{\#}$ by underlined Latin letters.

Endow the set of words of $F_d$ with the partial lexicographical
order. We put $x_{i_1}x_{i_2}\cdots x_{i_k}<x_{j_1}x_{j_2}\cdots
x_{j_t}$ if we have $i_1=j_1,\ldots,i_{s-1}=j_{s-1}$,
$i_{s}<j_{s}$ for some $s\geq1$. Note that if $v\in F_d^{\#}$,
then words $u$ and $uv$ are incomparable.

By an {\it identity} we mean an element of $K\LA F_d\RA$. All
identities are assumed to be $\mathcal{N}^{d}$-homogeneous, unless
the contrary is stated. The multidegree of an identity $t$ we
denote by $\mdeg(t)$. An identity $t$ is said to be an identity of
$N_{3,d}$, if the image of $t$ in $N_{3,d}$ under the natural
homomorphism is equal to zero. The zero polynomial is called the
{\it trivial} identity.

For identity $f=\sum_i\al_iu_i$, $\al_i\in K$, $u_i\in F_d^{\#}$,
$\ov{f}$ stands for the highest term of $f$, i.e., $\ov{f}$ is the
maximal word from the set $\{u_i\}$. It is easy to see that, due
to homogeneity, the highest term is unique. For a set of
identities $M$, denote by $\ov{M}$ the set of the highest terms of
the elements from $M$.

An identity is called {\it reduced} if the coefficient of its
highest term is equal to $1$.

We say that the identity is a {\it consequence} of a set of
identities, provided it belongs to the linear span of these
identities. As an example, we point out that $x_1^{3}$ is not a
consequence of $x_2^{3}$.

An element $\sum_i\al_i u_i\in K\LA F_d\RA$, where $\al_i\in K$,
$u_i\in F_d^{\#}$, is called {\it an element generated by words}
$v_1,\ldots,v_t$, if all $u_i$ are products of some elements from
$\{v_1,\ldots,v_t\}$.

For identities $t_1,t_2$ and for a set of identities $M$, notation
$t_1=t_2+\{M\}$ means that $t_1\in t_2+\lin M$.

Consider an element $g\in K\langle F_d\rangle$ and an identity
$t=u+(\sum_{i=1}^r \al_iu_i)$, where $\al_i\in K$,
$u,u_1,\ldots,u_r$ are pairwise different words. Let $g_1=g$. If
$g_k=v_1uv_2+\sum_{j=1}^s\be_j w_j$ for pairwise different words
$v_1uv_2,w_1,\ldots,w_s$ and $\be_j\in K$, then
$g_{k+1}=-\sum_{i=1}^r\al_iv_1u_iv_2+\sum_{j=1}^s\be_j w_j$. Note
that $g_{k+1}$ is not uniquely determined by $g_{k}$ and $t$. If
there is $k$ such that $g_k=\sum_{j=1}^s\be_j w_j$ for some words
$w_1,\ldots,w_s$ which do not contain subword $u$, then we say
that the chain $g_1,\ldots, g_k$ is finite and $g_k$ is its
result. If every chain is finite and they all have one and the
same result $g'$, then we call the identity $g'$ {\it the result
of application} of the identity $t$ with the {\it marked word} $u$
to the identity $g$. Otherwise we say that the result of
application of $t$ with the marked word $u$ to $g$ is indefinite.


The result of substitution $v_1\to u_1,\ldots,v_k\to u_k$ in $f\in
K\LA F_d \RA$, where $f$ is an element generated by words
$v_1,\ldots,v_k,v_{k+1},\ldots,v_{t}$, denote by $f|_{v_1\to
u_1,\ldots,v_k\to u_k}$. By {\it substitutional mapping} we mean
such homomorphism of $K$-algebras $\phi:K\LA F_k\RA\to K\LA
F_l\RA$ that $\phi(x_i)\in F_l^{\#}$, $i=\ov{1,k}$. A
substitutional mapping is called monotonous, if
$\phi(x_i)>\phi(x_j)$ for $x_i>x_j$. The set of monotonous
substitutional mappings denote by ${\mathcal M}_{k,l}$. Note that
for $\phi\in {\mathcal M}_{i,j}$, $\psi\in {\mathcal M}_{j,k}$ the
composition $\psi\circ\phi$ belongs to ${\mathcal M}_{i,k}$.

Denote $$T_1(a)=a^3,$$ $$T_2(a,b)=a^2b+aba+ba^2,$$
$$T_3(a,b,c)=abc+acb+bac+bca+cab+cba.$$ Partial and complete
linearization of the identity $f_1a^3f_2$ of $N_{3,d}$, where
$a\in K\LA F_d\RA^{\#}$, $f_1,f_2\in F_d$, gives that all
identities from $\sys=\{f_1T_1(a)f_2$, $f_1T_2(a,b)f_2$,
$f_1T_3(a,b,c)f_2|$ $a,b,c\in F_d^{\#}$, $f_1,f_2\in F_d\}$ are
identities of $N_{3,d}$. For  multidegree $\De$ let $\sys_{\De}$
be the subset of $\sys$ which consists of all identities of
multidegree $\De$. Clearly, each identity of $N_{3,d}(\De)$ is a
consequence of the set of identities $\sys_{\De}$. The set of
identities $\sys_{\De}$ can be treated like the system of
homogeneous linear equations in formal variables $\{w| w\in
F_d^{\#},\, \mdeg(w)=\De\}$. Then, free variables of the system
$\sys_{\De}$ form a basis for $\OK{\De}$. We call two systems of
linear equations  (two sets of identities, respectively)
equivalent, if the first one is a consequence of the second and
vice versa.

A word $w\in S$ is called {\it canonical} with respect to $x_i$,
if it has one of the following forms: $w_1$, $w_1x_iw_2$,
$w_1x_i^2w_2$, $w_1x_i^2ux_iw_2$, where words $w_1,w_2,u$ do not
contain $x_i$, words $w_1,w_2$ can be empty. If a word is
canonical with respect to all letters, then we call it {\it
canonical}. Number for future references the identity of $N_{3,d}$
\begin{eq}\label{eq_t2}
xux+(x^2u+ux^2), \; u\in F_d^{\#}.\;
\end{eq}

\section{Auxiliary results}


We will use the following facts from~\cite{Lopatin}:
\begin{lemma}\label{lemmaLop}
1. Applying identities~\Ref{eq_t2}, $x_iux_i^2=-x_i^2ux_i$  of
$N_{3,d}$, any non-zero word $w\in N_{3,d}$ can be represented as
a sum of canonical words which belong to the same homogeneous
component as $w$. In particular, if $\deg_{x_i}(w)>3$, $w\in F_d$,
then $w=0$ in $N_{3,d}$.

2. The inequality $x_1^2x_2^2x_1\neq0$ holds in $N_{3,d}$.

3. If $p=0$ or $p>3$, then $C(3,d,K)=6$ $(d\geq 2)$.

\quad If $p=2$, then $C(3,d,K)=d+3$, where $d\geq3$, and
$C(3,2,K)=6$.

4. If $p\neq3$, then $x^2ay^2=0$ is an identity of $N_{3,d}$,
where $a\in F_d^{\#}$.

5. If $p=3$, then $x^2y^2xay=x^2y^2xya$ is an identity of
$N_{3,d}$, where $a\in F_d^{\#}$.

6. If $p=2$, then  $x_1^2x_2\cdots x_dx_1\neq0$ holds in
$N_{3,d}$.

7. If $p\neq3$, then $I_1(x,a,b,c)=x^2abc+x^2acb$,
$I_2(x,a,b,c)=abcx^2+bacx^2$, $I_3(x,a,b,c)=ax^2bc+cax^2b$ are
identities of $N_{3,d}$, where $x$, $a$, $b$, $c$ are words.

\end{lemma}
\begin{proof}
1. See~\cite{Lopatin}, Statement 1.

2. See~\cite{Lopatin}, Statement 3.

3. See~\cite{Lopatin}, Propositions 1, 2.

4. See~\cite{Lopatin}, equality (5).

5. See~\cite{Lopatin}, proof of Statement 7.

6. See~\cite{Lopatin}, Statement 4.

7. Let $x,a,b,c$ be words, $p\neq3$. Partial linearization of the
identity from item $4$ with respect to $x$ ($y$, respectively)
gives that $I_1(x,a,b,c)$, $I_2(x,a,b,c)$ are identities of
$N_{3,d}$. Apply identity~\Ref{eq_t2}, where $x=x_1$, to the
identity $T_3(x_1a,bx_1,c)=0$ of $N_{3,d}$, and get that
$-T_3(x_1^2a,b,c)-T_3(a,bx_1^2,c)+3(bx_1^2ac+cbx_1^2a)=0$ in
$N_{3,d}$. Hence  $I_3(x,a,b,c)$ is the identity of $N_{3,d}$.
\end{proof}

\begin{remark}\label{remark3}
1. Consider a set $M=\{m_i\}_{i=\ov{1,s}}\subset K\LA
F_d\RA^{\#}$. Let $u\in F_d^{\#}$ be a word which is a summand of
one and only one element $m_1$ of the set $M$. Let $\sum_{i=1}^s
\al_im_i=0$ in $K\LA F_d\RA$, where $\al_i\in K$. Then $\al_1=0$.

2. Let $\De$ be a multidegree. Let $V=\{v_1,\ldots,v_s\}$ be a set
of words of  multidegree $\De$. Assume that for each word $w$, of
 multidegree $\De$, which do not lie in $V$ there is an
identity $w-f_w$ of $N_{3,d}$, where $f_w\in \lin V$. Then every
identity $\sum_i\al_iv_i$ ($\al_i\in K$) of $N_{3,d}(\De)$ is a
consequence of the identities which are results of application of
identities $\{w-f_w\}$ to the identities of $\sys_{\De}$. (Note
that results of these applications are defined.)
\end{remark}

\begin{lemma}\label{lemma1}
Let $d\geq1$. All identities, of $N_{3,d}(21^{d-1})$, generated by
$x_1^2,x_2,\ldots,x_d$ are consequences of the following
identities of $N_{3,d}(21^{d-1})$:

($a$) $f_1T_3(a,b,c)f_2$, where some word from $a$, $b$, $c$,
$f_1$, $f_2$ contains the subword $x_1^2$,

($b$) $3f_1I_i(x_1,a,b,c)f_2$, $i=\ov{1,3}$.

\noindent Here $a,b,c\in F_d^{\#}$, $f_1,f_2\in F_d$.
\end{lemma}
\begin{proof} By item~$2$ of Remark~\ref{remark3}, any identity, of
$N_{3,d}(21^{d-1})$, generated by $x_1^2,x_2,\ldots,x_d$ is a
consequence of identities which are results of application of
identity~\Ref{eq_t2} (where $x=x_1$) to the identities from
$\sys_{21^{d-1}}$.

If we apply~\Ref{eq_t2}, where $x=x_1$, to $T_2(x_1,a)=0$, then we
get a trivial identity.

The result of application of~\Ref{eq_t2}, where $x=x_1$, to an
identity $t=f_1T_3(a_1,a_2,a_3)f_2$, where $a_1,a_2,a_3\in
F_d^{\#}$, $f_1,f_2\in F_d$, denote by $t'$, and let
$t=\sum_{i=1}^6 u_i$ for some words $u_1,\ldots,u_6$. If words
$u_1,\ldots,u_6$ do not contain subword $x_1^2$ and $a_i\neq x_1$,
$i=\ov{1,3}$, then $t'$ is a consequence of identities ($a$). Let
$a,b,c\in F_d^{\#}$.

If $t=T_3(x_1,x_1,a)$, then $t'=0$.

If $t=T_3(x_1a,x_1,b)$, then $t'=-T_3(x_1^2,a,b)$.

If $t=T_3(ax_1b,x_1,c)$, then
$t'=-T_3(ab,x_1^2,c)-3I_3(x_1,a,b,c)$.

If $t=T_3(x_1a,bx_1,c)$, then
$t'=-T_3(x_1^2a,b,c)-T_3(a,bx_1^2,c)+3I_3(x_1,b,a,c)$.

If $t=x_1aT_3(x_1,b,c)$, then
$t'=-aT_3(x_1^2,b,c)-3I_1(x_1,a,b,c)$.

If $t=x_1T_3(x_1,a,b)$, then $t'=-T_3(x_1^2,a,b)$.

If $t=x_1T_3(x_1a,b,c)$, then
$t'=-x_1^2T_3(a,b,c)-T_3(x_1^2a,b,c)+3I_1(x_1,a,b,c)$.

Due to the fact that, if we read identities ($a$), ($b$) from
right to left, they do not change, the claim follows from the
regarded cases.
\end{proof}

Let $r=\ov{1,d}$. It is easy to see that for every $i=\ov{1,r}$
the result of application of the identity~\Ref{eq_t2}, where
$x=x_i$, to every identity  of multidegree $2^r1^{d-r}$ is
definite. For $\si\in S_r$ let $\psi_{\si}: K\LA
F_d\RA(2^r1^{d-r})\to K\LA F_d\RA(2^r1^{d-r})$ be such mapping
that $\psi_{\si}(t)$ is the result of the following procedure. Let
$t_1$ be the result of application of the identity~\Ref{eq_t2},
where $x=x_{\si(1)}$, to $t$. For $i=\ov{2,r}$ let $t_i$  be the
result of application of the identitiy~\Ref{eq_t2}, where
$x=x_{\si(i)}$, to $t_{i-1}$. We define $\psi_{\si}(t)=t_r$.

For any identity $t=\sum_{i=1}^s\al_iu_i$, $\al_i\in K$, $u_i\in
F_d^{\#}$, of multidegree $2^r1^{d-r}$, fix some permutations
$\si_{t,1},\ldots,\si_{t,s}\in S_r$. Consider the mapping $\psi:
K\LA F_d\RA(2^r1^{d-r})\to K\LA F_d\RA(2^r1^{d-r})$ such that for
$t=\sum_i\al_iu_i$, $\al_i\in K$, $u_i\in F_d^{\#}$, we have
$\psi(t)=\sum_i \al_i\psi_{\si_{t,i}}(u_i)$. Denote by $\Psi_r$
the set of all such mappings $\psi$.


\begin{lemma}\label{lemma16}
Let $d,r\geq1$, $\phi\in \Psi_r$. All identities, of
$N_{3,d}(2^{r}1^{d-r})$, generated by $x_1^2,\ldots,x_r^2$,
$x_{r+1}, \ldots, x_d$ are consequences of the following
identities of $N_{3,d}(2^r1^{d-r})$:

($a$) $f_1T_3(a,b,c)f_2$, where for each $k=\ov{1,r}$ some word
from $a$, $b$, $c$, $f_1$, $f_2$ contains the subword $x_k^2$,

($b$) $3\phi( f_1I_i(x_k,a,b,c)f_2)$, $i=\ov{1,3}$, $k=\ov{1,r}$,

($c$) the identity $3f_1x_i^2ax_j^2f_2$  ($i,j=\ov{1,r}$, $i\neq
j$) which is generated by $x_1^2,\ldots,x_r^2$, $x_{r+1}, \ldots,
x_d$.

\noindent Here $a,b,c\in F_d^{\#}$, $f_1,f_2\in F_d$.
\end{lemma}
\begin{proof}
For $\psi\in \Psi_r$ denote the sets of identities of multidegree
$2^r1^{d-r}$: $A_1=\{\psi_{\si}(w)-\psi_{\tau}(w)|\,w\in
F_d,\,\si,\tau\in S_r\}$,
$A_2^{\psi}=\{\psi(t)|\,t=f_1T_2(a,b)f_2,\, a,b\in F_d^{\#},\,
f_1,f_2\in F_d\}$ and
$A_3^{\psi}=\{\psi(t)|\,t=f_1T_3(a,b,c)f_2,\, a,b,c\in F_d^{\#},\,
f_1,f_2\in F_d\}$.

For an arbitrary $\psi\in \Psi_r$ identities, of
$N_{3,d}(2^{r}1^{d-r})$, generated by $x_1^2,\ldots,x_r^2$,
$x_{r+1}, \ldots, x_d$ are consequences of $A_1$, $A_2^{\psi}$,
$A_3^{\psi}$ (see item~$2$ of Remark~\ref{remark3}). The following
items conclude the proof.

1. {\it Identities $A_1$ are consequences of identities $(c)$. In
particular, $\psi(t)=\pi(t)+\{(a),(b),(c)\}$ and
$\psi(t+f)=\psi(t)+\psi(f)$ for any $\psi,\pi\in \Psi_r$, $t,f\in
K\LA F_d\RA(2^{r}1^{d-r})$}.

Proof. Note that the identity $\psi(3f_1x_i^2ax_j^2f_2)\in K\LA
F_d\RA(2^r1^{d-r})$, where $a\in F_d^{\#}$, follows from $(c)$.

If $r=1$ then $A_1=\{0\}$. Let $r=2$. Denote
$\psi_1(w)=\psi_{\si}(w)-\psi_{\tau}(w)$, where $\si=1\in S_2$,
$\tau=(1,2)\in S_2$. If $w=x_1ax_1bx_2cx_2$ or
$w=x_1ax_2bx_2cx_1$, $a,b,c\in F_d$, then the identity $\psi_1(w)$
is trivial. Consider $w=x_1ax_2bx_1cx_2$, $a,b,c\in F_d$. Then
$\psi_1(w)=(x_1ax_2bx_1)cx_2- x_1a(x_2bx_1cx_2)$ in $N_{3,d}$,
where the order of application of identity~\Ref{eq_t2} is
determined by parentheses.

If $ab,bc\neq 1$, then
$\psi_1(w)=(-x_1^2ax_2bcx_2-ax_2bx_1^2cx_2)-
(-x_1ax_2^2bx_1c-x_1abx_1cx_2^2)=
(x_1^2ax_2^2bc+x_1^2abcx_2^2+ax_2^2bx_1^2c+abx_1^2cx_2^2)-
(x_1^2ax_2^2bc+ax_2^2bx_1^2c+x_1^2abcx_2^2+abx_1^2cx_2^2)=0$.

If $a=b=c=1$, then $\psi_1(w)=0$.

If $a=b=1$, $c\neq 1$, then $\psi_1(w)=3x_1^2cx_2^2$.

If $b=c=1$, $a\neq 1$, then $\psi_1(w)=-3x_1^2ax_2^2$.

Therefore, if $r=2$ then the required is proved.

The case of $r>2$ follows from the case of $r=2$ and the fact that
any permutation is a composition of elementary transpositions.

2. {\it Identities $A_3^{\psi}$ are consequences of identities
$(a)$, $(b)$, $(c)$}.

Proof. It follows from Lemma~\ref{lemma1} and item~$1$.

3. {\it Identities $A_2^{\psi}$ are consequences of identities
$(a)$, $(b)$, $(c)$}.

Proof. We will use item~$1$ without reference. Prove by induction
on $k$ that for every identity $t=f_1T_2(v,a)f_2$ of multidegree
$2^r1^{d-r}$ we have $\psi(t)=0+\{(a),(b),(c)\}$, i.e.,
$\psi(f_1\cdot vav\cdot f_2)=-\psi(f_1\cdot v^2a\cdot
f_2)-\psi(f_1\cdot av^2\cdot f_2)+\{(a),(b),(c)\}$, where $v,a\in
F_d^{\#}$, $f_1,f_2\in F_d$, and $\deg(v)=k$.

Induction base is trivial.

Induction step. Without loss of generality we can assume that
$f_1,f_2$ are empty words. Consider a word $x_iu$ of  degree $k$,
where $i=\ov{1,r}$. We have
$\psi(T_2(x_iu,a))=\psi(x_iux_iua)+\psi(ax_iux_iu)+\psi(x_iuax_iu)+\{(a),(b),(c)\}$.
Induction hypothesis imply that
$\psi(x_iux_iua)=\psi(u^2x_i^2a)+\{(a),(b),(c)\}$,
$\psi(ax_iux_iu)=\psi(au^2x_i^2)+\{(a),(b),(c)\}$,
$\psi(x_iuax_iu)=\psi(x_i^2u^2a)+\psi(x_i^2au^2)+\psi(u^2ax_i^2)+\psi(ax_i^2u^2)+
\{(a),(b),(c)\}$. Thus,
$\psi(T_2(x_iu,a))=\psi(T_3(x_i^2,u^2,a))+\{(a),(b),(c)\}=0+\{(a),(b),(c)\}$
by item~$2$.
\end{proof}

\begin{lemma}\label{lemma12}
Let $p=3$. All identities of $\OK{1^d}$ are consequences of
identities $f_1T_3(a_1,a_2,a_3)f_2$ of multidegree $1^d$, where
$f_1,f_2\in F_d$, $a_1,a_2,a_3\in F_d^{\#}$, $\deg(a_1)\leq 3$,
$\deg(a_2)=\deg(a_3)=1$.
\end{lemma}
\begin{proof}
If $d\leq4$, then the statement is obvious.

Let $d\geq 5$. We prove by induction on $d$.

Induction base. In the case $d=5,6$ the statement was proven by
means of a computer programme.

Induction step. Let $d\geq 7$. Consider an identity
$t=a_1T_3(a_2,a_3,a_4)a_5$, $a_1,a_5\in F_d$, $a_2,a_3,a_4\in
F_d^{\#}$. There is such $k=\ov{1,5}$ that $\deg(a_k)\geq2$. Then
$a_k=x_ix_j\cdot w$ for some $w\in F_d$. Substituting $z$ for
subword $x_ix_j$, where $z$ is a new letter, and using induction
hypothesis, we get that $t\in\lin\{f_1T_3(b_1,b_2,b_3)f_2|$
$f_1,f_2\in F_d,\, b_1,b_2,b_3\in F_d^{\#},\, \deg(b_1b_2b_3)\leq
6\}$. The statement of the Lemma in the case $d=6$ concludes the
proof.
\end{proof}

A multilinear word $w=x_{\si(1)}\cdots x_{\si(d)}$, $\si\in S_d$,
is called {\it even}, if the permutation $\si$ is even, and {\it
odd} otherwise.

\begin{lemma}\label{lemma10}
$1$. Let $p=2,3$. Consider the homomorphism $\phi:K\LA
F_d\RA(1^d)\to K$, defined by $\phi(w)=1$, where $w\in F_d$. Then
$\phi$ maps identities of $N_{3,d}$ in zero.

$2$. Let $p=2$. For $i,j=\ov{1,d}$, $i\neq j$, consider the
homomorphism $\phi_{ij}:K\LA F_d\RA(1^d)\to K$, defined by the
following way: if $w=ux_ix_jv$, where $u,v\in F_d$, then
$\phi_{ij}(w)=1$, else $\phi_{ij}(w)=0$. Then $\phi_{ij}$ maps
identities of $N_{3,d}$ in zero.

$3$. Let $p=3$. Consider the homomorphism $\phi_{+}:K\LA
F_d\RA(1^d)\to K$, defined by the following way: for $w\in F_d$ we
have $\phi_{+}(w)=1$, if $w$ is even, else $\phi_{+}(w)=0$. Then
$\phi_{+}$ maps identities of $N_{3,d}$ in zero.

$4$. Let $p=3$. For $k=\ov{1,d}$ consider the homomorphism
$\phi_k: K\LA F_d\RA(\de_1,\ldots,\de_d)\to K\LA
F_d\RA(\de_1,\ldots,\de_{k-1},\de_{k+1},\ldots,\de_d)$, where
$\de_k=1,2$, defined by $\phi_k(x_{i_1}\cdots
x_{i_t})=y_{i_1}\cdots y_{i_t}$, where $y_{i_j}=x_{i_j}$, if
$i_j\neq k$, and $y_{i_j}=1$, if $i_j= k$ ($j=\ov{1,t}$,
$t=\delta_1+\cdots+\delta_d$). Then $\phi_k$ maps identities of
$N_{3,d}$ in zero.

$5$. Let $p=3$. For $k=\ov{1,d}$ consider the homomorphism \\
$\pi_k:K\LA
F_d\RA(\de_1,\ldots,\de_{k-1},3,\de_{k+1},\ldots,\de_d)\to K\LA
F_d\RA(\de_1,\ldots,\de_{k-1},1,\de_{k+1},\ldots,\de_d)$, defined
by the following way: for $w=u_1x_ku_2x_ku_3x_ku_4$, $u_i\in F_d$
$(i=\ov{1,4})$, we put
$\pi_k(w)=u_1(x_ku_2u_3+u_2x_ku_3+u_2u_3x_k)u_4$. Then $\pi_k$
maps identities of $N_{3,d}$ in zero, and
$\pi_{i}\pi_j(x_i^2x_j^2x_ix_j)=x_ix_j-x_jx_i$.
\end{lemma}
\begin{proof}
$2$. It is sufficient to proof that for
$t=b_1T_3(a_1,a_2,a_3)b_2$, where $b_1,b_2\in F_d$,
$a_1,a_2,a_3\in F_d^{\#}$,  we have $\phi_{ij}(t)=0$. If each word
from $\{b_1a_{\si(1)}a_{\si(2)}a_{\si(3)}b_2|\,\si\in S_3\}$ do
not contain subword $x_ix_j$, then $\phi_{ij}(t)=0$.

If there is $k$ such that $a_k=ux_ix_jv$, $u,v\in F_d$, then
$\phi_{ij}(t)=6=0$.

If $b_1=ux_i$, $a_1=x_jv$, $u,v\in F_d$, then $\phi_{ij}(t)=2=0$.

If $b_2=x_ju$, $a_1=vx_i$, $u,v\in F_d$, then $\phi_{ij}(t)=2=0$.

If $a_1=ux_i$, $a_2=x_jv$, $u,v\in F_d$, then $\phi_{ij}(t)=2=0$.

The statement follows from the regarded cases.

Items $3,4,5$ were proved in~\cite{Lopatin}; item $1$ is similar
to them.
\end{proof}

\section{The case of $p\neq3$, $d\leq3$}

\begin{prop}\label{prop1}
Let $p\neq3$, $d=\ov{1,3}$, $\De=(\de_1,\ldots,\de_d)$ is a
multidegree. Then, the set $B_{\De}$ is a basis for $\OK{\De}$,
where

\indent 1) the case of $\mod{\De}\leq3$:

\noindent $B_1=\{x_1\}$, $B_{1^2}=\{x_1x_2,\;x_2x_1\}$,
$B_{2}=\{x_1^2\}$,

\noindent $B_{1^3}=\{x_1x_2x_3,\; x_1x_3x_2,\; x_2x_1x_3,\;
x_2x_3x_1,\; x_3x_1x_2\}$, $B_{21}=\{x_1^2x_2,\; x_2x^2_1\;\}$;

\indent 2) the case of $\mod{\De}=4$:

\noindent $B_{21^2}=\{x_1^2x_2x_3,\; x_1^2x_3x_2,\; x_2x_1^2x_3,\;
x_2x_3x_1^2,\;x_3x_2x_1^2\}$,

\noindent $B_{2^2}=\{x^2_1x_2^2,\;x^2_2x_1^2\}$,
$B_{31}=\{x^2_1x_2x_1\}$;

\indent 3) the case of $\mod{\De}=5$:

\noindent $B_{2^21}=\{x^2_1x_2^2x_3,\; x^2_2x_1^2x_3,\;
x_3x_1^2x_2^2\}$,

\noindent $B_{31^2}=\{x^2_1x_2x_3x_1,\; x^2_1x_3x_2x_1\}$,

\noindent $B_{32}=\{x^2_1x_2^2x_1\}$;

\indent 4) otherwise $B_{\De}=\emptyset$.
\end{prop}
\begin{proof}
Cases of $\mod{\De}\leq3$ and $\De\in\{31,32\}$ follow from items
$1,2$ of Lemma~\ref{lemmaLop}.

If $\De\in\{21^2,2^2,2^21,31^2\}$, then we prove the statement by
considering the system of equations $\sys_{\De}$. Here we use item
$1$ of Lemma~\ref{lemmaLop}; and when $\De\in \{21^2,2^2,2^21\}$,
we use Lemma~\ref{lemma16} for decreasing the number of
considering equations.

Case $4)$ follows from item~$3$ of Lemma~\ref{lemmaLop}.
\end{proof}


\section{The case of $p=0$ or $p>3$}
We shall write $i$ for $x_i$, $i=\ov{1,d}$, so that it does not
lead to ambiguity.
\begin{prop}\label{prop_p0}
Let $p=0$ or $p>3$, $d\geq1$, $\De=(\de_1,\ldots,\de_d)$ is a
multidegree. Then the set $B_{\De}$ is a basis for $\OK{\De}$,
where

\indent 1) if $d\leq3$, then see Proposition~\ref{prop1},

\indent 2) $B_{1^4}=\left\{
\begin{array}{ccccc}
1234,& 1243,& 1324,& 1342,& 1423,\\
2134,& 2143,& 2314,& 2341,& 2413,\\
3124,& 3412&&&\\
\end{array}
\right\},$

\indent 3)  $B_{1^5}=\left\{
\begin{array}{ccccc}
12345,& 12354,& 12435,& 12453,& 12534,\\
13245,& 13254,& 13425,& 13452,& 13524,\\
14235,& 14523,& 23145,& 23415,& 23514
\end{array}
\right\},$

\indent 4) $B_{21^3}=\left\{
\begin{array}{cccc}
1^2234,& 1^2324,& 1^2423,&\\
21^234,& 21^243,& 231^24,& 241^23\\
\end{array}
\right\},$

\indent 5) if $\mod{\De}\geq6$, then $B_{\De}=\emptyset$.
\end{prop}
\begin{proof}
The  computations below were performed by means of a computer
programme. For $\De\in\{1^4,1^5,21^4\}$, we consider the
homogeneous system of linear equations $\sys_{\De}$ over the ring,
generated in ${\mathcal Q}$ by the set ${\mathcal
Z}\cup\{1/2,1/3\}$. Having expressed  higher words in terms of
lower words by the Gauss's method, we get that $\sys_{\De}$ is
equivalent to the system $\{1\cdot u=f_u|\, u\in
F_d,\,\mdeg(u)=\De,\,u\not\in B_{\De}\}$, where $f_u$ are linear
combinations of elements of $B_{\De}$. The statement is proven.

The case of $\mod{\De}\geq6$ follows from Lemma~\ref{lemmaLop}.
\end{proof}

\section{The composition method}

Denote by $M_5$ a basis for the space of identities of
$N_{3,5}(1^5)$ such that $M_5$ contains only reduced identities
and all elements of $M_5$ have the highest terms pairwise
different. The basis of this kind exists, because if  $t_1,t_2$
are reduced identities with $\ov{t}_1=\ov{t}_2$ and $t_1\neq t_2$,
then $\ov{t}_1\neq\ov{t_1-t_2}$ and
$\lin\{t_1,t_2\}=\lin\{t_1,t_1-t_2\}$. For $d\geq5$, let
$$M_d=\{t\in K\LA F_d\RA(1^d)|\mbox{ there are } t'\in M_5,\,\phi\in {\mathcal M_{5,d}},\,
a\in F_d\mbox{ such that } t=a\phi(t')\}.$$ %
Identities from $M_d$ are identities of $N_{3,d}$.

Let $\De=(\de_1,\ldots,\de_d)$ be a multidegree. For a set
$J\subset K\LA F_{d}\RA(\De)$, denote $B(J)=\{w\in F_d^{\#}|\,
\mdeg(w)=\De, w\not\in \ov{J}\}$. Since every word, of multidegree
$\De$, which do not belong to $B(J)$, can be expressed in terms of
lower words by applying identities of $J$; therefore, for any
$f\in K\LA F_d\RA(\De)$, we have
\begin{eq}\label{eq4}
f=\sum\nolimits_i\al_i t_i+\sum\nolimits_j\be_j w_j, \mbox{ where
} \al_i,\be_j\in K,\, t_i\in J,\, w_j\in B(J),\, \ov{t_i},w_j\leq
\ov{f}.
\end{eq}
Thus,
\begin{eq}\label{eq5}
\frac{K\LA F_d\RA(\De)}{\lin J}\simeq\lin B(J), \mbox{ and, in
particular, }\OK{1^d}=\lin \Phi(B(M_d)),
\end{eq}

\noindent where $\Phi:K\LA F_d\RA^{\#}\to N_{3,d}$ is the natural
homomorphism. Further, we will write $B(M_d)$ instead of
$\Phi(B(M_d))$ so that it does not lead to ambiguity.

\noindent{\bf Definition.} A set of reduced identities $M$ is
called  {\it complete under composition}, if for any $t_1,t_2\in
M$, where $\ov{t}_1=\ov{t}_2$, we have $t_1-t_2=\sum_{i=1}^k \al_i
g_i$, $\al_i\in K$, $g_i\in M$, and $\ov{g}_i<\ov{t}_1$,
$i=\ov{1,k}$.

\begin{lemma}\label{lemma_id}
For $d\geq5$, all identities of $\OK{1^d}$ are consequences of
identities of $M_d$.
\end{lemma}
\begin{proof}
For any identity $t$ of $\OK{1^d}$, we have $t=\sum_i\al_i t_i$
for some $\al_i\in K$, $t_i=f_iT_3(a_i,b_i,c_i)g_i$, $f_i,g_i\in
F_d$, $a_i,b_i,c_i\in F_d^{\#}$. For any $i$, exists $\phi_i\in
{\mathcal M}_{5,d}$ and an identity $t_i'$, of $N_{3,5}(1^5)$,
such that $t_i=\phi_i(t_i')$. The set $M_5$ is a basis for the
space of identities of $N_{3,5}(1^5)$, thus
$t_i'=\sum_j\be_{ij}g_{ij}$, where $g_{ij}\in M_5$, $\be_{ij}\in
K$. Hence $\phi_i(t_i')=\sum_j\be_{ij}\phi_i(g_{ij})\in \lin M_d$.
Therefore, $t\in\lin M_d$.
\end{proof}
\begin{lemma}\label{lemma_compos}
(Composition Lemma~\cite{Bokut}) For $d\geq5$ $B(M_d)$ is a basis
for $\OK{1^d}$ if and only if $M_d$ is complete under composition.
\end{lemma}
\begin{proof}
Let $B(M_d)$ be a basis for $\OK{1^d}$. For $t_1,t_2\in M_d$,
where $\ov{t}_1=\ov{t}_2$, let $g=t_1-t_2$. By formula~\Ref{eq4}
we have $g=\sum_i\al_if_i+\sum_j\be_jw_j$, where
$\ov{f}_i\leq\ov{g}<\ov{t}_1$, $f_i\in M_d$, $w_i\in B(M_d)$,
$\al_i,\be_j\in K$. The identity $g-\sum_i \al_if_i$ is an
identity of $N_{3,d}$, $B(M_d)$ is a basis for $\OK{1^d}$, hence
$\sum_j\be_jw_j=0$ in $K\LA F_d\RA$. Thus $g=\sum_i\al_if_i$, and
the claim is proven.

Let $M_d$ be complete under composition. Assume that, on the
contrary, $B(M_d)$ is not a basis. Thus, the set $B(M_d)$ is
linearly dependent in $N_{3,d}$ (see equality~\Ref{eq5}). Hence
there is a non-trivial identity $f=\sum_i\al_iu_i$, $\al_i\in K$,
$u_i\in B(M_d)$, such that $f=0$ in $N_{3,d}$. Note that
$\ov{f}\not\in \ov{M}_d$.

Lemma~\ref{lemma_id} implies $f=\sum_{j=1}^k \be_j t_j$, where
$\be_j\in K^{\ast}$, $t_j\in M_d$. Without loss of generality, we
can assume that for some $s$, we have
$\ov{t}_1=\cdots=\ov{t}_s=\ov{f}$. If $s=1$ then we get a
contradiction to $\ov{f}\not\in \ov{M}_d$. Let $s\geq2$. Since
$M_d$ is complete under composition, for $j=\ov{2,s}$ we have
$t_1-t_j=\sum_l \ga_{jl} g_{jl}$, where $g_{jl}\in M_d$,
$\ov{g}_{jl}<\ov{t}_1$, $\ga_{jl}\in K$. Expressing $t_j$ from
these equalities, we get $f=\la t_1+\sum_q\la_q h_q$ for some
$h_q\in M_d$, $\ov{h}_q<\ov{t}_1$, $\la,\la_q\in K$. If
$\la\neq0$, then $\ov{f}\in \ov{M}_d$, so we get a contradiction.
Thus, $\la=0$. Repeating the same argument several times, we get a
contradiction to the non-triviality of $f$.
\end{proof}

\begin{lemma}\label{lemma3}
Let $d\geq 5$. Then for any $t_1,t_2\in M_d$, where
$\ov{t}_1=\ov{t}_2$, there are $s=\ov{5,10}$, $t'_1,t'_2\in M_s$,
$\phi\in {\mathcal M}_{s,d}$, $a\in F_d$ such that
$t_i=a\phi(t'_i)$, where  $i=1,2$.
\end{lemma}
\begin{proof}
By definition of $M_d$, there are $t_i''\in M_5$, $\psi_i\in
{\mathcal M}_{5,d}$, $c_i\in F_d$ such that
$t_i=c_i\psi_i(t_i'')$, where $i=1,2$. Denote
$w=\ov{t}_1=\ov{t}_2$. Let $\ov{t}_1''=x_{j_1}\cdots x_{j_5}$,
$\ov{t}_2''=x_{k_1}\cdots x_{k_5}$. Consider a partition of $w$
into subwords $w=a\cdot a_1\cdot\,\cdots\, \cdot a_s$, where
$s=\ov{5,10}$, $a_1,\ldots,a_s\in F_d^{\#}$, $a\in F_d$, which is
the result of intersection of partitions $w=\ov{t}_1=c_1\cdot
\psi_1(x_{j_1})\cdot\,\cdots\,\cdot \psi_1(x_{j_5})$,
$w=\ov{t}_2=c_2\cdot \psi_2(x_{k_1})\cdot\,\cdots\,\cdot
\psi_2(x_{k_5})$. Here we assume $a\neq1$ if and only if $c_1,c_2$
are non-empty words. Then $c_i=ad_i$, where $d_i\in F_d$, $i=1,2$.
There is a permutation $\si\in S_s$ such that if $x_i<x_j$, then
$a_{\si(i)}<a_{\si(j)}$, where $i,j=\ov{1,s}$. Since
$d_1\psi_1(t_1''),d_2\psi_2(t_2'')$ are elements generated by
words $a_1,\ldots,a_s$, the substitutions
$d_i\psi_i(t_i'')|_{a_{\si(j)}\to x_j,j=\ov{1,s}}=t_i'$, $i=1,2$,
are well-defined. It is easy to see that for $t_i'$ there is
$\psi_i'\in {\mathcal M}_{5,s}$, $e_i\in F_s$ such that
$e_i\psi_i'(t_i'')=t_i'$, i.e. $t_i'\in M_s$, where $i=1,2$.
Define $\phi\in{\mathcal M}_{s,d}$ by the following way:
$\phi(x_j)=a_{\si(j)}$, $j=\ov{1,s}$. We have $a\phi(t_i')=t_i$,
where $i=1,2$, thus the claim is proven.
\end{proof}

\begin{prop}\label{prop_lemma4}
If $B(M_d)$ is a basis for $\OK{1^d}$ for $d=\ov{5,10}$, then
$B(M_d)$ is a basis for $\OK{1^d}$ for any $d\geq 5$.
\end{prop}
\begin{proof}
Let $d\geq 11$. Consider $t_1,t_2\in M_d$, where
$\ov{t}_1=\ov{t}_2$. We apply Lemma~\ref{lemma3} to $t_1,t_2$;
further we use the notation from Lemma~\ref{lemma3}. By the data
and Lemma~\ref{lemma_compos} $M_s$ is complete under composition,
hence $t'_1-t'_2=\sum_i \al_ig_i'$, where $\al_i\in K$, $g_i'\in
M_s$, $\ov{t}_1'>\ov{g}_i'$. Thus $t_1-t_2=\sum_i\al_ig_i$, where
$g_i=a\phi(g_i')\in M_d$; because of the composition of monotonous
substitutional mappings is a monotonous substitutional mapping. By
monotony of $\phi$, we have $\ov{t}_1>\ov{g}_i$ for all $i$. Hence
$M_d$ is complete under composition, and Lemma~\ref{lemma_compos}
concludes the proof.
\end{proof}

\section{Multilinear homogeneous component}\label{section6}

\noindent{\bf Notation.} For $p=2,3$ and $d\geq1$ recursively
define sets $B_{1^d}$ of the words of multidegree $1^d$.

Let $p=2$. Then

$\begin{array}{l}
1) \;B'_1=\{x_1\};\\
2) \mbox{ for } d\geq2\;\mbox{ define }
B'_{1^d}=x_1\{B_{1^{d-1}}|_{x_i\to x_{i+1},\,i=\ov{1,d-1}}\}
\cup \{\un{e}_{d,k}|\,k=\ov{2,d}\} \cup \{\un{f}_{d}\} \cup \\
\{\un{h}_{d,k}|\,k=\ov{3,d}\}.\\
\end{array}$

Here, if $i\geq 1,i\neq4$, then $B_{1^i}=B_{1^i}'$;
$B_{1^4}=B_{1^4}'\cup \{x_2x_1x_4x_3\}$ and

$\begin{array}{l} \un{e}_{d,k}=x_2\cdots x_k\cdot x_1\cdot
x_{k+1}\cdots x_d\; (d\geq2,\,
k=\ov{2,d});\\
\un{f}_d=x_2\cdots x_{d-2}\cdot x_{d}x_1x_{d-1}\; (d\geq3);\\
\un{h}_{d,k}=x_k\cdot x_1\cdots \widetilde{x}_{k}\cdots x_d\;
(d\geq3,\,
k=\ov{3,d}).\\
\end{array}$

Let $p=3$. Then

$\begin{array}{l}
1)\; B_1=\{x_1\};\\
2)  \mbox{ for } d\geq2\;\mbox{ define }
B_{1^d}=x_1\{B_{1^{d-1}}|_{x_i\to x_{i+1},\,i=\ov{1,d-1}}\}\cup
x_2\{B_{1^{d-1}}|_{x_i\to x_{i+1},\, i=\ov{2,d-1}}\}\cup\\
\{\un{e}_{d,k}|\,k=\ov{3,d}\}.
\end{array}$

Here, $\un{e}_{d,k}=x_3\cdots x_k\cdot x_1x_2\cdot x_{k+1}\cdots
x_d\; (d\geq3,\, k=\ov{3,d})$.

For future needs define $B_{1^0}=\{1\}$, where $1$ stands for the
empty word.

The aim of this section is to prove the following theorem:
\begin{theo}\label{theo_Bd}
For $p=2,3$, $d\geq1$ the set $B_{1^d}$ is a basis for $\OK{1^d}$.
\end{theo}

\begin{remark} It is not difficult to see that%
$$|B_{1^d}|=\left\{
\begin{array}{lcc}
d(d-1)&,& p=2,\,d\geq4\\
2^d-d&,& p=3\\
\end{array}\right..$$
\end{remark}

Let $V$ be a finite dimensional vector space over $K$,
$V=\lin\{v_1,\ldots,v_m\}$, where non-zero vectors $\{v_i\}$ are
linearly ordered by the following way: $v_1<\cdots<v_m$. Note that
the vectors $v_1,\ldots,v_m$ need not be linearly independent.

\noindent {\bf Definition.} A basis of $V$
$v_{k_1},\ldots,v_{k_s}$ is called {\it minimal} (with respect to
the linearly ordered set $v_1,\ldots,v_m$), if for any
$i=\ov{1,m}$ we have $v_i=\sum_j \al_{ij}v_{k_j}$, $\al_{ij}\in
K$, where $k_j\leq i$.

Consider  $L_1=\{v_{j_1},\ldots,v_{j_s}\}$,
$L_2=\{v_{k_1},\ldots,v_{k_s}\}$ which are bases for $V$, where
$j_1<\cdots<j_s$, $k_1<\cdots<k_s$. We write $L_1<L_2$, if there
is $l=\ov{1,s}$ such that
$j_1=k_1,\ldots,j_{l-1}=k_{l-1},j_{l}<k_{l}$.

\begin{lemma}\label{lemma5}
$1$. A basis of $V$ $v_{k_1},\ldots,v_{k_s}$ is minimal (with
respect to the linearly ordered set $v_1,\ldots,v_m$) if and only
if it is the least one with respect to the determined linear
order.

$2$. The minimal basis is uniquely determined.
\end{lemma}
\begin{proof}
$1$. Let $L$ be the minimal basis. Then
\begin{eq}\label{eq6}
\mbox{ if }v_k\not\in\lin\{v_1,\ldots,v_{k-1}\}, \mbox{ then }
v_k\in L; \mbox{ else } v_k\not\in L  \; (k=\ov{1,m}).
\end{eq}
Thus, $L$ is the least basis.

Let $L\subset\{v_1,\ldots,v_m\}$ be the least basis. Then
condition~\Ref{eq6} is valid for it. Thus $L$ is the minimal
basis.

$2$. This item follows from item $1$.
\end{proof}

\indent Apply aforesaid on minimal bases to $\OK{1^d}$. As a
linearly ordered set we take $\{w\in F_d|\, \mdeg(w)=1^d\}$.
\begin{lemma}\label{lemma6}
$1$. $B(M_5)$ is a basis for $\OK{1^5}$.

$2$. Let $d\geq5$. If $B(M_d)$ is a basis for $\OK{1^d}$,  then
$B(M_d)$ is the minimal basis.
\end{lemma}
\begin{proof}
$1$. The definition of $M_5$ implies that $M_5$ is complete under
composition. Lemma~\ref{lemma_compos} concludes the proof.

$2$. Identities $M_d$ imply that any word, of $N_{3,d}$, which do
not belong to $B(M_d)$ can be expressed in terms of lower words.
Thus $B(M_d)$ is the minimal basis.
\end{proof}

\begin{lemma}\label{lemma7}
$1$. If $p=2$, then for $i=4,5$ the minimal basis for
$N_{3,i}(1^i)$ is $B_{1^i}$, and

$\begin{array}{rcl} B_{1^4}&=&\left\{
\begin{array}{ccccc}
1234,&1243,&1324,&1342,&1423,\\
2134,&2143,&2314,&2341,&2413,\\
3124,&4123\\
\end{array}
\right\};\\
B_{1^5}=B(M_5)&=&\left\{
\begin{array}{ccccc}
12345,&12354,&12435,&12453,&12534,\\
13245,&13254,&13425,&13452,&13524,\\
14235,&15234,\\
21345,&23145,&23415,&23451,\\
23514,\\
31245,&41235,&51234\\
\end{array}
\right\}.
\end{array}$

$2$. If $p=3$, then for $i=4,5$ the minimal basis for
$N_{3,i}(1^i)$ is $B_{1^i}$, and

$\begin{array}{rcl} B_{1^4}&=&\left\{
\begin{array}{ccccc}
1234,&1243,&1324,&1342,&1423,\\
2134,&2143,&2314,&2341,&2413,\\
3124,&3412\\
\end{array}
\right\};\\
B_{1^5}=B(M_5)&=&\left\{
\begin{array}{ccccc}
12345,&12354,&12435,&12453,&12534,\\
13245,&13254,&13425,&13452,&13524,\\
14235,&14523,\\
21345,&21354,&21435,&21453,&21534,\\
23145,&23154,&23415,&23451,&23514,\\
24135,&24513,\\
31245,&34125,&34512\\
\end{array}
\right\}.
\end{array}$
\end{lemma}
\begin{proof}
By Lemma~\ref{lemma6}, $B(M_5)$ is the minimal basis for
$N_{3,5}(1^5)$, and, in particular, it does not depend on the
choice of $M_5$ (see Lemma~\ref{lemma5}). Expressing higher words
in terms of lower words by Gauss's method, we solve the system
$\sys_{1^d}$ and find the minimal basis for $\OK{1^d}$. These
calculations were performed by means of a computer programme for
$d=4,5$, $p=2,3$.
\end{proof}

\begin{lemma}\label{lemma8}
If $p=2$, $d\geq5$, then

\noindent $\ov{M}_d= \left\{\begin{array}{ll}
w_3=ua_1a_2a_3,& a_1>a_2>a_3,\\
w_{4,1}=ua_1a_2a_3a_4,& a_1>a_2,a_4\mbox{ and } a_3>a_4,\\
w_{4,2}=ua_1a_2a_3a_4,& a_1>a_4\mbox{ and } a_2>a_3,\\
w_{5,1}=ua_1a_2a_3a_4a_5,& a_1>a_2\,\;\mbox{ and }\;\; a_3>a_4\mbox{ or }a_3>a_5 \mbox{ or }a_4>a_5,\\
w_{5,2}=ua_1a_2a_3a_4a_5,& a_1>a_3\,\;\mbox{ and }\;\; a_2>a_4\mbox{ or } a_2>a_5 \mbox{ or }a_4>a_5,\\
w_{5,3}=ua_1a_2a_3a_4a_5,& a_1>a_4\mbox{ and } a_2>a_5.\\
\end{array}
\right\},$

\noindent where all elements of $\ov{M}_d$ are words of
multidegree $1^d$, $u\in F_d$, $a_i\in F_d^{\#}$, $i=\ov{1,5}$.

If $p=3$, $d\geq5$, then

\noindent $\ov{M}_d= \left\{\begin{array}{ll}
w_3=ua_1a_2a_3,& a_1>a_2>a_3,\\
w_{4}=ua_1a_2a_3a_4,& a_1>a_2,a_4, \\
w_{5,1}=ua_1a_2a_3a_4a_5,& a_1>a_2,a_3\mbox{ and }a_4>a_5,\\
w_{5,2}=ua_1a_2a_3a_4a_5,& a_1>a_3,a_4\mbox{ and }a_2>a_5,\\
w_{5,3}=ua_1a_2a_3a_4a_5,& a_1>a_4,a_5\mbox{ and }a_2>a_3.\\
\end{array}
\right\},$

\noindent where all elements of $\ov{M}_d$ are words of
multidegree $1^d$, $u\in F_d$, $a_i\in F_d^{\#}$, $i=\ov{1,5}$.
\end{lemma}
\begin{proof} It is sufficient to prove the statement for $d=5$.
Denote by $M$ the set from the formulation of the Lemma.
Considering all possibilities, we get that $B(M_5)=\{w\in
F_5|\,\mdeg(w)=1^5,\,w\not\in M\}$ (see Lemma~\ref{lemma7}).
\end{proof}

\begin{lemma}\label{lemma9}
For $p=2,3$, $d\geq5$ we have $B(M_d)=B_{1^d}$.
\end{lemma}
\begin{proof}
For $p=2,3$ we get $B_{1^d}\subset B(M_d)$ by Lemma~\ref{lemma8}.

Further $w_{i,j},w_i$ stand for the words from Lemma~\ref{lemma8}.

{\bf The case of  ${\mathbf {p=2}}$}. Inclusion $B(M_d)\subset
B_{1^d}$ follows from items $1$, $2$ (see below) by induction on
$d$.

1. {\it If $w=x_{i_1}\cdots x_{i_d}\in B(M_d)$, $i_1\geq 3$, then
$w=\un{h}_{d,i_1}$.}

Proof. Let $w=x_{i_1}u$. The word $w$ contains letters $x_1,x_2$.
Denote by $u_1,u_2,u_3$ some elements of $F_d$. If
$u=u_1x_2u_2x_1u_3$, then $w=w_3\in \ov{M}_d$, that is a
contradiction. Hence $u=u_1x_1u_2x_2u_3$. If the word $u_1$ is not
empty, then $w=w_{4,2}\in\ov{M}_d$, that is a contradiction. If
the word $u_2$ is not empty, then $w=w_{4,1}\in\ov{M}_d$, that is
a contradiction. Assume that $u_3=x_{j_1}\cdots x_{j_s}$ and there
are $k,t$ such that $k<t\leq s$ and $j_k>j_t$. Then
$w=w_{5,1}\in\ov{M}_d$, that is a contradiction. Hence
$w=\un{h}_{d,i_1}$.

2. {\it If $w=x_2x_{i_2}\cdots x_{i_d}\in B(M_d)$, then
$w=\un{e}_{d,k}$ for some $k=\ov{2,d}$ or $w=\un{f}_d$.}

Proof. Consider $w=x_2 u_1 x_1 u_2$, where $u_1,u_2\in F_d$. If
there are not $r,s$ such that $r>s$ and $u_1u_2=v_1x_rv_2x_sv_3$,
where $v_1,v_2,v_3\in F_d$, then $w=\un{e}_{d,k}$ for some
$k=\ov{2,d}$. Assume that there are such $r,s$. If the word $v_3$
contains the letter $x_1$, then $w=w_{4,2}\in\ov{M}_d$; a
contradiction. If the word $v_1$ contains the letter $x_1$, then
$w=w_{5,1}$ or $w=w_{5,2}$, hence $w\in\ov{M}_d$; a contradiction.
Let the word $v_2$ contains the letter $x_1$. If the word $v_3$ is
not empty, then $w=w_{5,2}\in\ov{M}_d$; a contradiction. If
$\deg(v_2)>1$, then $w=w_{5,2}$ or $w=w_{5,3}$, hence
$w\in\ov{M}_d$; a contradiction. There is the only possibility
which we have not considered, namely $w=\un{f}_d$.


{\bf The case of $\mathbf {p=3}$}. Inclusion $B(M_d)\subset
B_{1^d}$ follows from items $1$, $2$ (see below) by induction on
$d$.

1. {\it If $w=x_{i_1}\cdots x_{i_d}$, $i_1\geq 4$, then $w\not\in
B(M_d)$.}

Proof. The word $w_1=x_{i_2}\cdots x_{i_d}$ contains the letters
$x_1,x_2,x_3$. There are $r,s=\ov{1,3}$ such that the word $w_1$
contains some letters between letters $x_r,x_s$. We also have
$x_{i_1}>x_r,x_s$. Thus $w=w_{4}\in\ov{M}_d$, that is $w\not\in
B(M_d)$.

2. {\it If $w=x_3x_{i_2}\cdots x_{i_d}\in B(M_d)$, then
$w=\un{e}_{d,k}$  for some $k=\ov{3,d}$.}

Proof. If $w=x_3u_1x_2u_2x_1u_3$ for some elements $u_1,u_2,u_3$
of $F_d$, then $w=w_3\in \ov{M}_d$; a contradiction. Thus
$w=x_3u_1x_1u_2x_2u_3$. If the word $u_2$ is not empty, then
$w=w_{4}\in\ov{M}_d$, that is a contradiction. Hence
$u=x_3u_1x_1x_2u_3$. If there are $r,s$ such that $r>s$ and
$u_1u_3=v_1x_rv_2x_sv_3$ for some $v_1,v_2,v_3\in F_d$, then
$w=w_{5,1}$ or $w=w_{5,2}$ or $w=w_{5,3}$. Therefore $w\in
\ov{M}_d$, that is a contradiction. Hence there are not such
$r,s$. So $w=\un{e}_{d,k}$ for some $k=\ov{3,d}$.
\end{proof}

\begin{lemma}\label{lemma11}
Let $p=3$, $d\geq6$, $\phi_k$ is the mapping from
Lemma~\ref{lemma10}, where $k=\ov{1,d}$. Then for any $w\in
B(M_{d})$ we have $\phi_k(w)|_{x_i\to x_{i-1},i=\ov{k+1,d}}\in
B(M_{d-1})$.
\end{lemma}
\begin{proof}
If for a word $w$ of  multidegree $1^d$ we have
$\phi_k(w)|_{x_i\to x_{i-1},i=\ov{k,d}}\in \ov{M}_{d-1}$, then
$w\in \ov{M}_d$, that is $w\not\in B(M_d)$.
\end{proof}

\begin{proof_theo1}
If $d=1,2,3$ then obviously $B_{1^d}$ is a basis. For $d=4,5$
$B_{1^d}$ is a basis by Lemma~\ref{lemma7}.
Proposition~\ref{prop_lemma4} and Lemma~\ref{lemma9} imply that in
order to prove the Theorem it is sufficient to verify that
$B_{1^d}$ is linearly independent in $N_{3,d}$ for $d=\ov{6,10}$.
This verification was done by means of a computer programme
applying the algorithm described below.

{\bf The case of $\mathbf{p=2}$}. Assume that there is an identity
$f=\sum_{w\in B_{1^d}}\al_w w$, $\al_w\in K$, such that $f=0$ in
$N_{3,d}$. Considering $\phi_{ij}(t)=0$, $\phi(t)=0$, where
$\phi_{ij},\phi$ are mappings from Lemma~\ref{lemma10}, we get a
homogeneous system of linear equations in $\{a_w\}$. Having solved
this system we get that $\al_w=0$ for any $w\in B_{1^d}$. It was
calculated  for $d=\ov{6,10}$ by means of a computer programme.

{\bf The case of $\mathbf {p=3}$}. By Lemma~\ref{lemma9}, we have
$B_{1^d}=B(M_d)$. Identities from $M_5$ express elements of the
set $\{w\in F_5|\, \mdeg(w)=1^5\}$ in terms of elements of
$B_{1^5}$. Applying these identities, we get that for any word $w$
of multidegree $1^d$ $w=f_w$ in $N_{3,d}$, where $f_w\in\lin
B_{1^d}$. Applying identities  $\{w=f_w\}$, rewrite identities
$\{f_1T_3(a_1,a_2,a_3)f_2|\,f_1,f_2\in F_d,\, \deg(a_1)\leq 3,
\deg(a_2)=\deg(a_3)=1\}$ in terms of linear combinations of the
elements of $B_{1^d}$. As a result, we get only trivial
identities. This, together with Lemma~\ref{lemma12}, imply that
the system of identities $\sys_{1^d}$ is equivalent to the set of
identities $M=\{w-f_w|\,w\in F_d,\mdeg(w)=1^d,\,w\not\in
B_{1^d}\}$. The set $M$ is linearly independent in $K\LA F_d\RA$,
thus $B_{1^d}$ is linearly independent in $N_{3,d}$. The given
algorithm performed by means of a computer programme proved that
$B_{1^d}$ is linearly independent in $N_{3,d}$ when $d=\ov{6,9}$.
Here we need Lemma~\ref{lemma12} in order to decrease the quantity
of identities which have to be considered. For $d=10$ described
algorithm ran for a long time, thus we used another approach to
the case of $d=10$.

Assume that there is an identity $f=\sum_{w\in B_{1^d}}\al_w w$,
where $\al_w\in K$, such that $f=0$ in $N_{3,d}$. Consider
mappings $\phi_{k}$ $(k=\ov{1,d})$, $\phi_{+}$ from
Lemma~\ref{lemma10}. We assume that $B_{1^{d-1}}$ is linearly
independent in $N_{3,d}$. Hence, we get a homogeneous system of
linear equations in $\{\al_w\}$ (see also Lemma~\ref{lemma11}).
For even $d=\ov{6,10}$ it was calculated by means of a computer
programme that this system has the only solution $\al_w=0$ for any
$w\in B_{1^d}$.
\end{proof_theo1}

\section{The case of $p=2$}
Let $d\geq4$, $i,j=\ov{2,d}$, $i\neq j$. Introduce notations for
some words of multidegree $21^{d-1}$: $\un{a}_i=x_1^2x_ix_2\cdots
\widetilde{x}_{i}\cdots x_d$, $\un{b}_i=x_1\cdots
\widetilde{x}_{i}\cdots x_dx_ix_1^2$,
$\un{c}_{ij}=x_ix_1^2x_jx_2\cdots \widetilde{x}_{i}\cdots
\widetilde{x}_{j}\cdots x_d$, if $i<j$, and
$\un{c}_{ij}=x_ix_1^2x_jx_2\cdots \widetilde{x}_{j}\cdots
\widetilde{x}_{i}\cdots x_d$, if $i>j$. By items~$1$,~$7$ of
Lemma~\ref{lemmaLop} we have
\begin{eq}\label{eq0}
\OK{21^{d-1}}=\lin\{\un{a}_i,\un{b}_i,\un{c}_{ij}|\,i,j=\ov{2,d},\,i\neq
j\}.
\end{eq}
\begin{theo}\label{theo2}
Let $p=2$, $d\geq1$.

1. For $\De=(\de_1,\ldots,\de_d)$, $d\leq3$ a basis for $\OK{\De}$
is the set $B_{\De}$ defined in Proposition~\ref{prop1}.

2. For $d\geq4$ a basis for $\OK{1^d}$ is the set $B_{1^d}$
defined above.

3. A basis for $N_{3,4}(21^3)$ is the set
$B_{21^{3}}=\{\un{a}_i,\un{b}_i,\un{c}_{23},\un{c}_{32}|\,
i=\ov{2,4}\}$.

\quad For $d\geq5$ a basis for $\OK{21^{d-1}}$ is the set
$B_{21^{d-1}}=\{\un{a}_i,\un{b}_i,\un{c}_{23}|\, i=\ov{2,d}\}$.

4. For $d\geq4$ a basis for $\OK{2^21^{d-2}}$ is the set
$B_{2^21^{d-2}}=\{x_1^2x_2^2x_3\cdots x_d, x_2^2x_1^2x_3\cdots
x_d\}$.

5. For $d\geq4$ a basis for $\OK{31^{d-1}}$ is the set
$B_{31^{d-1}}=\{x_1^2x_2\cdots x_dx_1\}$.

6. The rest of ${\mathcal N}^d$-homogeneous components of
$N_{3,d}$ are equal to zero.
\end{theo}

In order to prove item $3$ we need the following Lemma.

Denote $h_{ij}=\un{b}_i+\un{c}_{ij}+\un{a}_j$. Consider identities
of multidegree $21^{d-1}$:
$$\begin{array}{rcl}
M_0&=&\{f_1I_i(x_1,a,b,c)f_2 |\,i=\ov{1,3}\},\\
M_1&=&\mbox{the set of identities
(a) from Lemma~\ref{lemma1}},\\
M_2&=&\{f_1T_3(x_1^2,a,b)f_2\},\\
M_3&=&\{T_3(x_1^2,x_i,x_j)a, \; aT_3(x_1^2,x_i,x_j) |\, 2\leq i<j\leq d\},\\
M_4&=&\left\{
\begin{array}{l}
\{h_{23}+h_{34},h_{23}+h_{42},h_{32}+h_{24},h_{32}+h_{43}\}, \mbox{if\, } d=4\\
\{h_{23}+h_{ij}|\, 2\leq i\neq j\leq d,\, i\neq2\mbox{ or }j\neq3\}, \mbox{if\, } d\geq5,\\
\end{array}
\right. \\
\end{array}$$
where $a,b,c\in F_d^{\#}$, $f_1,f_2\in F_d$. For $i=\ov{1,4}$
denote $L_i=\lin\{M_0\cup M_i\}$.

\begin{lemma}\label{lemma2}
Let $p=2$, $d\geq4$. Then

1. $L_1=L_2$.

2. $L_2=L_3$.

3. $L_3=L_4$.
\end{lemma}
\begin{proof}
Let $i,j,k,l\in \ov{2,d}$ be pairwise different numbers.

1. Inclusion $L_2\subset L_1$ is obvious.

Consider an identity $t\in M_1$.

If $t=T_3(ax_1^2b,c,d)$, then
$t=aI_1(x_1,b,c,d)+I_2(x_1,c,d,a)b+I_3(x_1,ca,b,d)+dI_3(x_1,a,b,c)+2dcax_1^2b$.

If $t=T_3(x_1^2a,b,c)$, then
$t=I_1(x_1,a,b,c)+I_3(x_1,c,a,b)+I_3(x_1,b,a,c)$.

If $t=x_1^2aT_3(b,c,d)$, then
$t=I_1(x_1,a,b,c)d+I_1(x_1,a,b,d)c+I_1(x_1,a,c,d)b$.

If $t=x_1^2T_3(a,b,c)$, then
$t=I_1(x_1,a,b,c)+I_1(x_1,b,a,c)+I_1(x_1,c,a,b)$.

Thus we get $L_1\subset L_2$.

3.  Introduce notations for identities:
$$\begin{array}{l}
f_{ijk}=\un{a}_i+\un{a}_j+\un{c}_{ij}+\un{c}_{ji}+\un{c}_{ik}+
\un{c}_{jk}=T_3(x_1^2,x_i,x_j)x_ka+\{M_0\}, \\
g_{ijk}=\un{b}_j+\un{b}_k+\un{c}_{jk}+\un{c}_{kj}+\un{c}_{ij}+
\un{c}_{ik}=ax_iT_3(x_1^2,x_j,x_k)+\{M_0\},\\
\end{array}$$
where $a\in F_d$.  We have
\begin{eq}\label{eq1}
f_{ijk}+g_{jik}=h_{ki}+h_{ij}\in L_3.
\end{eq}
\indent  If $d=4$, then formula~\Ref{eq1} implies $L_4\subset
L_3$.

If $d\geq5$, then replacing indices in formula~\Ref{eq1} we get
$h_{ij}+h_{jk}\in L_3$, $h_{ij}+h_{jl}\in L_3$, $h_{jl}+h_{lk}\in
L_3$, $h_{lk}+h_{kj}\in L_3$. Add up last four formulas and get
$h_{jk}+h_{kj}\in L_3$. Thus $h_{ij}+h_{ji},h_{ij}+h_{jk}\in L_3$.
Therefore $L_4\subset L_3$.

From $f_{ijk}=4h_{23}+h_{ij}+h_{ji}+h_{ik}+h_{jk}\in L_4$,
$g_{ijk}=4h_{23}+h_{jk}+h_{kj}+h_{ik}+h_{ij}\in L_4$ we can see
that $L_3\subset L_4$ for $d\geq5$.

Let $d=4$. Equalities $f_{ijk}=f_{jik}$, $g_{ijk}=g_{ikj}$,
$f_{234}+f_{243}=f_{342}$, $g_{234}+g_{324}=g_{423}$,
$$\begin{array}{ccl}
f_{234}&=&(h_{23}+h_{34})+(h_{32}+h_{24}),\\
f_{243}&=&(h_{23}+h_{42})+(h_{32}+h_{24})+(h_{32}+h_{43}),\\
g_{234}&=&(h_{23}+h_{34})+(h_{32}+h_{24})+(h_{32}+h_{43}),\\
g_{324}&=&(h_{23}+h_{34})+(h_{23}+h_{42})+(h_{32}+h_{24}),\\
\end{array}$$
imply $L_3\subset L_4$.

2. Inclusion $L_3\subset L_2$ is obvious.

Consider an identity $t\in M_2$.

If $t=T_3(x_1^2,a,bc)$, then
$t=T_3(x_1^2,a,b)c+bT_3(x_1^2,a,c)+\{M_0\}$.

If $d\geq5$, then for $t=f_1x_kT_3(x_1^2,x_i,x_j)x_lf_2$,
$f_1,f_2\in F_d$, we have
$t=\un{c}_{ki}+\un{c}_{kj}+\un{c}_{ij}+\un{c}_{ji}+\un{c}_{il}+\un{c}_{jl}+
\{M_0\}\in L_4$; thus $t\in L_3$ by item~$3$.

Therefore $L_2\subset L_3$.
\end{proof}

\begin{proof_theo2}
{2.} See Theorem~\ref{theo_Bd}.

{ 3.} By equality~\Ref{eq0} and identities $M_4$ it is sufficient
to show that $B_{21^{d-1}}$ is linear independent in $N_{3,d}$.

Lemmas~\ref{lemma1} and~\ref{lemma2} imply that all identities, of
$\OK{21^{d-1}}$, generated by $x_1^2,x_2,\ldots,x_d$ are
consequences of identities $M_0\cup M_4$. Let $f$ be a mapping
such that the image of a word $w$, generated by
$x_1^2,x_2,\ldots,x_d$, of multidegree $21^{d-1}$ is equal to the
result of application of the identities from item~$7$ of
Lemma~\ref{lemmaLop} to $w$, i.e. $f(w)$ is equal to $\un{a}_i$,
$\un{b}_i$ or $\un{c}_{ij}$ for some $i,j$. Identities $M_0\cup
M_4$ are equivalent to the identities $M_0'\cup M_4=M$, where
$M_0'=\{w+f(w)| w$ is a word, $\mdeg(w)=21^{d-1}$, $w\neq f(w)\}$.
Every identity of $M_0'$ ($M_4$, respectively) contains a word
which is not a summand of any element of $B_{21^{d-1}}$ and is a
summand of one and only one identity of $M$ ($M_4$, respectively).
Moreover we can assume that regarded words are pairwise different.
Thus item~$2$ of Remark~\ref{remark3} implies that $B_{21^{d-1}}$
is linearly independent in $N_{3,d}$.


{ 4.} Denote $a=x_1^2x_2^2x_3\cdots x_d$, $b=x_2^2x_1^2x_3\cdots
x_d$. By items~$1$,~$4$ and~$7$ of Lemma~\ref{lemmaLop}, we have
$\lin\{a,b\}=\OK{2^21^{d-2}}$. We claim that $a,b$ are linearly
independent in $N_{3,d}$. Consider the homomorphism of vector
spaces $\psi: K\LA F_d\RA\to K\LA F_d\RA(2^21^{d-2})$ defined by
the following way: for a word $w$
 $\psi(w)=\al a+\be b$, where $\al$ ($\be$, respectively)
is equal to the number of subwords  $x_1x_2$ ($x_2x_1$,
respectively) in the word $w$.

For $u,v\in F_d^{\#}$ define%
$$
\psi(u,v)=\left\{
\begin{array}{ccl}
a&,& \mbox{if } u=u_1x_1,\, v=x_2v_1\,(u_1,v_1\in F_d)\\
b&,& \mbox{if } u=u_1x_2,\, v=x_1v_1\,(u_1,v_1\in F_d)\\
0&,& \mbox{otherwise }
\end{array}\right..$$%
It is easy to see that for $u_1,\ldots,u_s\in F_d^{\#}$
\begin{eq}\label{eq9}
\psi(u_1\cdots
u_s)=\sum\limits_{i=1}^s\psi(u_i)+\sum_{i=1}^{s-1}\psi(u_i,u_{i+1}).
\end{eq}

Consider an identity $t$ of $\sys_{2^21^{d-2}}$.

If $t=f_1T_1(g)f_2$, $f_1,f_2\in F_d$, $g\in F_d^{\#}$, then
$t\not\in\sys_{2^21^{d-2}}$. It is a contradiction.

If $t=f_1T_2(g_1,g_2)f_2$, $f_1,f_2\in F_d$, $g_1,g_2\in
F_d^{\#}$, then
$\psi(t)=\psi(g_2)+\psi(f_1)+\psi(f_2)+\psi(f_1,g_2)+\psi(g_2,f_2)$,
by equality~\Ref{eq9}. The multidegree of $t$ is $2^21^{d-2}$,
thus $g_1\in\{x_1,x_2,x_1x_2,x_2x_1\}$. Hence $\psi(t)=0$.

If $t=f_1T_3(g_1,g_2,g_3)f_2$, $f_1,f_2\in F_d$, $g_1,g_2,g_3\in
F_d^{\#}$, then $\psi(t)=0$ by equality~\Ref{eq9}.

Therefore $\psi(t)=0$ for any identity $t$ of $\OK{2^21^{d-2}}$.

If $t=\al a+\be b$, $\al,\be\in K$, is an identity of $N_{3,d}$,
then $\psi(t)=\al a+\be b=0$ in $K\LA F_d\RA$. Hence $\al=\be=0$.

{ 5.} By Lemma~\ref{lemmaLop} we have $u=x_1^2x_2\cdots x_dx_1\neq
0$. Identities from items~$1$,~$7$ of Lemma~\ref{lemmaLop} imply
that $x_1^2abcx_1=x_1(abcx_1^2)=x_1bacx_1^2=x_1^2bacx_1$ in
$N_{3,d}$. The last identity together with item~$7$ of
Lemma~\ref{lemmaLop} imply that $\lin\{u\}=N_{3,d}(31^{d-1})$.

{6.} It follows from item~$3$ of Lemma~\ref{lemmaLop}.
\end{proof_theo2}

\section{The case of $p=3$}

{\bf Notation.} For $p=3$, $r,s,l\geq0$ determine the set
$B_{3^r1^s}$ of words of multidegree $3^r1^s$ and the set
$B_{3^r2^s1^l}$ of words of  multidegree $3^r2^s1^l$:

$B_{3^{2r}1^s}=\un{u}_{2r}\{B_{1^s}|_{x_i\to x_{i+2r},\,
i=\ov{1,s}}\}\cup\{\un{q}_{2r,s,k}|\,k=\ov{1,s}\}$ ($r\geq 0$),

$B_{3^{2r+1}1^s}=\un{u}_{2r}x_{2r+1}^2x_{2r+2}\{B_{1^s}|_{x_1\to
x_{2r+1},\,x_i\to x_{i+2r+1},\,
i=\ov{2,s}}\}\cup\{\un{q}_{2r+1,s,k}|\,k=\ov{3,s+1}\}\cup
\{\un{q}_{2r+1,s}\}$ ($r\geq0$), where

$\un{q}_{2r,s,k}=\un{u}_{2r-2}\cdot x_{2r-1}^2x_{2r}^2\cdot
x_{2r+1}\cdots x_{2r+k} \cdot x_{2r-1}x_{2r}\cdot x_{2r+k+1}\cdots
x_{2r+s}$   ($r,s\geq1$, $k=\ov{1,s}$),

$\un{q}_{2r+1,s,k}=\un{u}_{2r}\cdot x_{2r+1}^2\cdot x_{2r+3}\cdots
x_{2r+k}\cdot x_{2r+1}x_{2r+2}\cdot x_{2r+k+1} \cdots x_{2r+s+1}$
($r\geq0$, $s\geq2$, $k=\ov{3,s+1}$),

$\un{q}_{2r+1,s}=\un{u}_{2r-2}\cdot
x_{2r-1}^2x_{2r}^2x_{2r+1}^2x_{2r-1}x_{2r}x_{2r+1}\cdot
x_{2r+2}\cdots x_{2r+s+1}$ ($r\geq1$, $s\geq0$),

$\un{u}_{2k}=\un{w}_{12}\cdots\un{w}_{2k-1,2k}$ $(k\geq1)$,
$\un{u}_0$ is the empty word,

$\un{w}_{ij}=x_i^2x_j^2x_ix_j$.

Define $B_{3^r2^s1^l}=B_{3^r1^{s+l}}|_{x_i\to
x_{i}^2,\,i=\ov{r+1,r+s}}$.

As an example we point out that $B_{3^01^s}=B_{1^s}$,
$B_{3^{2r}}=\{\un{u}_{2r}\}$, $B_{3^{2r+1}}=\{\un{q}_{2r+1,0}\}$,
$B_{3^{2r+1}1}=\{\un{u}_{2r}x_{2r+1}^2x_{2r+2}x_{2r+1},\,
\un{q}_{2r+1,1}\}$.

\begin{theo}\label{theo3}
Let $p=3$.

$1$. A basis for $N_{3,d}(3^r2^s1^l)$ is the set $B_{3^r2^s1^l}$,
where $r,s,l\geq0$.

$2$. The rest of ${\mathcal N}^d$-homogeneous components of
$N_{3,d}$ are equal to zero.
\end{theo}

\begin{remark} For $p=3$, $r\geq2$, $s,l\geq0$ we have $|B_{3^r2^s1^l}|=2^{s+l}$.
\end{remark}

\begin{proof_theo3}
1. Consider the homomorphism $\phi:\OK{3^{r}1^{s+l}}\to
\OK{3^r2^s1^l}$,
defined by%
$$\phi(x_i)=\left\{
\begin{array}{ccl}
x_i^2&,& r+1\leq i\leq r+s\\
x_i&,& {\rm otherwise}\\
\end{array}
\right..
$$ %
By item~$1$ of Lemma~\ref{lemmaLop} $\phi$ is surjective. By
Lemma~\ref{lemma16} $\phi$ is injective. Thus, $\phi$ is an
isomorphism of vector spaces. Hence it is sufficient to prove the
Theorem for multidegree $3^r1^s$. The last follows from
Lemmas~\ref{lemma15},~\ref{lemma14} (see below).

2. See item~$1$ of Lemma~\ref{lemmaLop}.
\end{proof_theo3}

Further for multilinear elements $f_1,f_2\in K\langle F_d\rangle$,
where $\deg(f_1f_2)=m$,  writing $f_1\xi(f_2)$ means that $\xi$ is
a substitutional mapping of ${\mathcal M}_{d,m}$ such that the
multidegree of $f_1\xi(f_2)$ is equal to $1^m$.

Consider $i_1,\ldots,i_r\in\{1,2\}$, where $r<d$, and a word $u$
such that
$$%
\begin{array}{lcl}
u=x_1&,& \mbox{if } r=d-1\\
u\in\{\un{e}_{d-r,k}|\,k=\ov{3,d-r}\}&,& \mbox{if } r<d-1\\
\end{array}.$$%
Denote by $(i_1i_2\ldots i_r;u)$ the word $w\in B_{1^d}$ which is
the result of the following procedure. Let $w_{r+1}=u\in
B_{1^{d-r}}$, $w_{k}=x_{i_{k}}\xi(w_{k+1})\in B_{1^{d-k+1}}$ for
every $k=\ov{1,r}$. Put $w=w_1$. For short, we will write
$(1^{s-1}i_s\ldots i_r;u)$ instead of $(1\ldots1 i_s\ldots i_r;u)$
and so on.

\begin{lemma}\label{lemma15}
Let $r,s\geq0$. Then $\lin B_{3^{r}1^s}=\OK{3^r1^s}$.
\end{lemma}
\begin{proof}
Let $J_{r,s}=M_{2r+s} \cup \{ ux_{2i-1}x_{2i}v, \;
ux_{2i}x_{2i-1}v, \; u(x_{2i}ax_{2i-1}+x_{2i-1}ax_{2i})v, \;\\
u(x_{2i-1}x_{2j-1}x_{2i}bx_{2j} - x_{2i-1}x_{2j-1}x_{2i}x_{2j}b)v
\,|\, u,v\in F_{2r+s},\; a,b\in F_{2r+s}^{\#},\; i,j=\ov{1,r}, \;
b>x_{2j} \}$ be the subset of $K\LA F_{2r+s}\RA(1^{2r+s})$.

Consider the homomorphism $\phi:K\LA F_{2r+s}\RA(1^{2r+s})\to
\OK{3^r1^s}$, defined by $\phi(x_{2i-1})=x_i^2$,
$\phi(x_{2i})=x_i$, $\phi(x_{j})=x_{j-r}$, where $i=\ov{1,r}$,
$j=\ov{2r+1,2r+s}$. By item~$1$ of Lemma~\ref{lemmaLop}, $\phi$ is
surjective. Identities  $x^2y^2xay=x^2y^2xya$, $xax^2+x^2ax=0$ of
$N_{3,d}$ (see Lemma~\ref{lemmaLop}) imply that $\phi$ induces the
epimorphism $\phi_1:K\LA F_{2r+s}\RA(1^{2r+s})/\lin(J_{r,s})\to
\OK{3^r1^s}$. By equality~\Ref{eq5} we have
$\OK{3^r1^s}=\lin\,\phi_1(B(J_{r,s}))$. So in order to prove the
statement it is sufficient to prove that
\begin{eq}\label{eq7}
\phi_1(B(J_{r,s}))=B_{3^r1^s}.
\end{eq}
\indent Note that $B(J_{r,s})=B(J_{r-1,s+2})\setminus\ov{J}_{r,s}$
for $r\geq1$, $s\geq0$.

For $r=0$  equality~\Ref{eq7} is obvious.

Let $r=1$. We have $B(J_{1,s})=B_{1^{s+2}}\setminus \ov{J}_{1,s}=
x_1 \xi(B_{1^{s+1}}) \cup x_2 \xi(B_{1^{s+1}})\cup \{x_3\cdots
x_k\cdot x_1x_2\cdot x_{k+1}\cdots x_{s+2}|\,
k=\ov{3,s+2}\}\setminus \ov{J}_{1,s}= x_1 \xi(B_{1^{s+1}})
\setminus \ov{J}_{1,s}= x_1x_2 \xi(B_{1^s})\cup
x_1x_3\xi(B_{1^s})\cup \{x_1\cdot x_4\cdots x_k\cdot x_2x_3\cdot
x_{k+1}\cdots x_{s+2}|\,k=\ov{4,s+2}\} \setminus \ov{J}_{1,s}=
x_1x_3\xi(B_{1^s})\cup \{x_1\cdot x_4\cdots x_k\cdot x_2x_3\cdot
x_{k+1}\cdots x_{s+2}|\,k=\ov{4,s+2}\}$. Hence equality~\Ref{eq7}
holds.

Let $r=2$. We have $B(J_{2,s})=B(J_{1,s+2})\setminus
\ov{J}_{2,s}=/\mbox{see above}/= x_1x_3\xi(B_{1^{s+2}})\cup
\{x_1\cdot x_4\cdots x_k\cdot x_2x_3\cdot x_{k+1}\cdots
x_{s+4}|\,k=\ov{4,s+4}\}\setminus \ov{J}_{2,s}=
x_1x_3\xi(B_{1^{s+2}}) \setminus \ov{J}_{2,s}=
x_1x_3x_2\xi(B_{1^{s+1}}) \cup x_1x_3x_4\xi(B_{1^{s+1}})\cup
\{x_1x_3\cdot x_5\cdots x_k\cdot x_2x_4\cdot x_{k+1}\cdots
x_{s+4}|\,k=\ov{5,s+4}\} \setminus \ov{J}_{2,s}=
x_1x_3x_2x_4\xi(B_{1^s})\cup x_1x_3x_2x_5\xi(B_{1^s})\cup
\{x_1x_3x_2\cdot x_6\cdots x_k\cdot x_4x_5\cdot x_{k+1}\cdots
x_{s+4}|\,k=\ov{6,s+4}\} \cup \{x_1x_3\cdot x_5\cdots x_k\cdot
x_2x_4\cdot x_{k+1}\cdots x_{s+4}|\, k=\ov{5,s+4}\}
\setminus\ov{J}_{2,s} =x_1x_3x_2x_4 \xi(B_{1^s}) \cup
\{x_1x_3\cdot x_5\cdots x_k\cdot x_2x_4\cdot x_{k+1}\cdots
x_{s+4}|\,k=\ov{5,s+4}\}$. Thus  equality~\Ref{eq7} holds.

Let $r=3$. We have $B(J_{3,s})=B(J_{2,s+2})\setminus
\ov{J}_{3,s}=/\mbox{see above}/= x_1x_3x_2x_4\xi(B_{1^{s+2}})\\
\cup \{x_1x_3\cdot x_5\cdots x_{k}\cdot x_2x_4\cdot x_{k+1}\cdots
x_{s+6}|\,k=\ov{5,s+6}\}
 \setminus \ov{J}_{3,s}=
x_1x_3x_2x_4\xi(B(J_{1,s})) \cup \{x_1x_3\cdot x_5\cdot
x_2x_4\cdot x_6\cdots x_{s+6}\}$. Thus equality~\Ref{eq7} is
proved to be true.

Let $r\geq4$. By induction on $r$ prove that
\begin{eq}\label{eq8}
B(J_{r,s})=x_1x_3x_2x_4\xi(B(J_{r-2,s})), \mbox{ where } r\geq4.
\end{eq}
\indent Induction base. Let $r=4$. We have $B(J_{4,s})=
B(J_{3,s+2})\setminus \ov{J}_{4,s}=/\mbox{see above}/=
x_1x_3x_2x_4\xi(B(J_{1,s+2}))\cup \{x_1x_3x_5x_2x_4x_6\cdot
x_7\cdots x_{s+8}\}\setminus \ov{J}_{4,s}=
x_1x_3x_2x_4\xi(B(J_{2,s}))$.

Induction step. We have $B(J_{r,s})= B(J_{r-1,s+2})\setminus
\ov{J}_{r,s}= /$induction hypothesis$/=
x_1x_3x_2x_4\xi(B(J_{r-3,s+2}))\setminus \ov{J}_{r,s}=
x_1x_3x_2x_4\xi(B(J_{r-2,s}))$.

Formula~\Ref{eq8} implies that equality~\Ref{eq7} is valid.
\end{proof}

\begin{lemma}\label{lemma14}
Let $r,s\geq0$. Then

1. The set $B_{3^{2r}1^s}$ is linearly independent in $N_{3,d}$,
where $d=2r+s$.

2. The set $B_{3^{2r+1}1^s}$  is linearly independent in
$N_{3,d}$, where $d=2r+s+1$.
\end{lemma}
\begin{proof}
Denote by $\pi_i$ the homomorphism from Lemma~\ref{lemma10}. Note
that for any $k\geq1$, $v\in B_{1^k}$ words $x_1x_2\xi(v)$,
$x_2x_1\xi(v)$ belong to $B_{1^{k+2}}$. Denote $\pi_1\cdots
\pi_{2r-2}(\un{u}_{2r-2})=u$.

1. If $r=0$, then see Theorem~\ref{theo_Bd}. Let $r\geq1$. For
$v\in B_{1^s}$, $k=\ov{1,s}$ we have
$\pi_{1}\cdots\pi_{2r}(\un{u}_{2r}\xi(v))=u
(x_{2r-1}x_{2r}-x_{2r}x_{2r-1})\xi(v)$,
$\pi_{1}\cdots\pi_{2r}(\un{q}_{2r,s,k})=u\xi(a-b-c+d)$, where
$a=x_3\cdots x_{k+2}\cdot x_1x_2\cdot x_{k+3}\cdots x_{s+2} =
\un{e}_{s+2,k+2}\in B_{1^{s+2}}$, $b=x_1\cdot x_3\cdots
x_{k+2}\cdot x_2 \cdot x_{k+3}\cdots x_{s+2}=(12^k1^{s-k};x_1)\in
B_{1^{s+2}}$, $c=x_2\cdot x_3\cdots x_{k+2}\cdot x_1 \cdot
x_{k+3}\cdots x_{s+2}=(2^{k+1}1^{s-k};x_1)\in B_{1^{s+2}}$,
$d=x_1x_2\cdot x_3\cdots x_{k+2}\cdot x_{k+3}\cdots
x_{s+2}=(1^{s+1};x_1)\in B_{1^{s+2}}$. By above remark, for any
$w\in B_{3^{2r}1^s}$ we have $\pi_1\cdots
\pi_{2r}(w)=\sum_i\al_{w,i}a_{w,i}$, where $\al_{w,i}\in K$,
$a_{w,i}\in B_{1^{2r+s}}$. By item~$1$ of Remark~\ref{remark3}, we
get that the set $\{\sum_i\al_{w,i}a_{w,i}|\,w\in B_{3^{2r}1^s}
\}$ is linearly independent in $K\LA F_d\RA$. So the assumption
that $B_{3^{2r}1^s}$ is linearly dependent in $N_{3,d}$ gives that
$B_{1^{2r+s}}$ is linearly dependent in $N_{3,d}$ (see
Lemma~\ref{lemma10}), and the last contradicts
Theorem~\ref{theo_Bd}.

2. Let $v\in B_{1^s}$, $k=\ov{3,s+1}$. Denote $a_v=\pi_1\cdots
\pi_{2r+1}( \un{u}_{2r} x_{2r+1}^2 x_{2r+2} \xi(v))$,
$b_k=\pi_1\cdots \pi_{2r+1}(\un{q}_{2r+1,s,k})$, $c=\pi_1\cdots
\pi_{2r+1}(\un{q}_{2r+1,s})$. Let $\phi_1$ be the homomorphism
from item~$4$ of Lemma~\ref{lemma10}. We have that
$a_v=u(x_{2r-1}x_{2r}-x_{2r}x_{2r-1})(x_{2r+2}\xi(v)-x_{2r+1}x_{2r+2}\xi(\phi_1(v)))$,
$b_k=u(x_{2r-1}x_{2r}-x_{2r}x_{2r-1})(x_{2r+3}\cdots x_{2r+k}\cdot
x_{2r+1}x_{2r+2} \cdot x_{2r+k+1}\cdots x_{2r+s+1}-x_{2r+1}\cdot
x_{2r+3}\cdots x_{2r+k}\cdot x_{2r+2}\cdot x_{2r+k+1}\cdots
x_{2r+s+1})=%
u(x_{2r-1}x_{2r}-x_{2r}x_{2r-1})\xi(\un{e}_{s+1,k}-(12^{k-2}1^{s-k+1};x_1))$,
$c=u(x_{2r-1}x_{2r+1}x_{2r}-x_{2r}x_{2r-1}x_{2r+1}+x_{2r}x_{2r+1}x_{2r-1}-
x_{2r+1}x_{2r-1}x_{2r})x_{2r+2}\cdots x_{2r+s+1}=
u\xi((121^{s};x_1)-(21^{s+1};x_1)+(2^21^s;x_1)-\un{e}_{s+3,3})$.
Above remark and Lemma~\ref{lemma11} imply that $a_v,b_k,c\in \lin
B_{1^{2r+s+1}}$. Thus $a_v=\sum_i\al_{v,i}a_{v,i}$,
$b_k=\sum_i\be_{k,i}b_{k,i}$, $c=\sum_i \ga_ic_i$, where
$\al_{v,i},\be_{k,i},\ga_i\in K$, $a_{v,i},b_{k,i},c_i\in
B_{1^{2r+s+1}}$. The set $\{a_v,b_k,c|\,v\in
B_{1^s},\,k\in\ov{3,s+1}\}$ is linearly independent, because of
each set from  the class of sets
$\{\{a_{v,i}\},\{b_{k,i}\},\{c_i\}|\,v\in
B_{1^s},\,k\in\ov{3,s+1}\}$ contains an element which does not
belong to other sets. Thus the assumption that $B_{3^{2r+1}1^s}$
is linearly dependent in $N_{3,d}$ gives that $B_{1^{2r+s+1}}$ is
linearly dependent in $N_{3,d}$ (see Lemma~\ref{lemma10}), and the
last contradicts Theorem~\ref{theo_Bd}.
\end{proof}

\section{Matrix invariants}

Let $n\geq2$. Denote by $M_{n,d}(K)=M_{n}(K) \oplus\cdots \oplus
M_{n}(K)$ the sum of $d$ copies of the space of $n\times n$
matrices. The general linear group $GL_n(K)$ acts on $M_{n,d}(K)$
by diagonal conjugation: for $g\in GL_n(K)$, $A_i\in M_n(K)$
$(i=\ov{1,d})$ we have $g(A_1,\ldots,A_d) = (gA_1g^{-1}, \ldots,
gA_dg^{-1})$. The coordinate ring of the affine space $M_{n,d}(K)$
is the polynomial algebra $K_{n,d}=K[x_{ij}(r)|\,1\leq i,j\leq
n,\; r=\ov{1,d}]$, where $x_{ij}(r)$ stands for the function such
that the image of $(A_1,\ldots,A_d)\in M_{n,d}(K)$ is $(i,j)$th
entry of the matrix $A_r$. The action of $GL_n(K)$ on $M_{n,d}(K)$
induces the action on $K_{n,d}$: $(g\cdot f)(A)=f(g^{-1}A)$, where
$g\in GL_n(K)$, $f\in K_{n,d}$, $A\in M_{n,d}$. Denote by
$R_{n,d}=\{f\in K_{n,d}|\mbox{ for all } g\in GL_n(K): gf=f\}$ the
matrix algebra of invariants. Let $X_r=(x_{ij}(r))_{1\leq i,j\leq
n}$ be the generic matrices of order $n$ $(r=\ov{1,d})$, and let
$\si_k(A)$ be the coefficients of the characteristic polynomial of
a matrix $A\in M_n(K)$, that is $\det(\la
E-A)=\la^n-\si_1(A)\la^{n-1}+\cdots+(-1)^n\si_n(A)$. The algebra
$R_{n,d}$ is generated by all elements of the form
$\si_k(X_{i_1}\cdots X_{i_s})$ (see~\cite{Don}). The
Procesi--Razmyslov Theorem on the relations in $R_{n,d}$ was
extended to the case of a field of an arbitrary characteristic
in~\cite{Zub96}.

The goal of the constructive theory of invariants is to find a
minimal (i.e. irreducible) homogeneous system of generators
(shortly m.h.s.g.) for the algebra of invariants.  A m.h.s.g. for
$R_{2,d}$ was determined in~\cite{Sibirskii} for $p=0$,
in~\cite{Procesi} for $p>2$, and in~\cite{DKZ} for $p=2$.
In~\cite{Dom} some upper and lower bounds on the highest degree of
elements of a m.h.s.g. for $R_{n,d}$ are pointed out for an
arbitrary $p$. In~\cite{Abeasis} in the case $p=0$ the cardinality
of a m.h.s.g. for $R_{3,d}$ was calculated  for $d\leq 10$ on a
computer, and was shown a way how such set can be constructed by
means of a computer programme. The explicit upper bound on the
highest degree of elements of a m.h.s.g. for $R_{3,d}$ is given
in~\cite{Lopatin} (except for the case $p=3$, $d=6k+1$, $k>0$,
where the least upper bound is estimated with error not greater
than $1$). In this section we point out a m.h.s.g. for $R_{3,d}$
for an arbitrary $p$, $d$.

The algebra $R_{n,d}$ possesses natural ${\mathcal N}$- and
${\mathcal N}^d$-gradings by degrees and multidegrees
respectively.  Denote by $R_{n,d}^{+}$ the subalgebra generated by
all elements of $R_{n,d}$ of positive degree. An element $r\in
R_{n,d}$ is called {\it decomposable}, if it can be expressed in
terms of elements of $R_{n,d}$ of lower degree, that is it belongs
to the ideal $(R_{n,d}^{+})^2$. Clearly, $\{r_i\}\in R_{n,d}$ is a
m.h.s.g. if and only if $\{\ov{r}_i\}$ is a basis for
$\ov{R}_{n,d} =R_{n,d} / (R_{n,d}^{+})^2$. If two elements
$r_1,r_2\in R_{n,d}$ are equal modulo the ideal $(R_{n,d}^{+})^2$,
we write $r_1\equiv r_2$. There is a close connection between
decomposability of an element of $R_{n,d}$ and equality to zero of
some element of $N_{n,d}$ (see Lemma~\ref{lemmaLop2} below).  Let
$A_{n,d}$ be a $K$-algebra without unity, generated by the generic
matrices $X_1,\ldots,X_d$. The homomorphism of algebras
$\Phi:A_{n,d}\to N_{n,d}$, defined by $\Phi(X_i)=x_i$,  is defined
correctly (see~\cite{Lopatin}).

\begin{remark}\label{remark1}
For each $\De=(\de_1,\ldots,\de_d)$, where $\de_1\geq \cdots
\geq\de_d\geq0$, let $\MSG_{\De}\subset R_{n,d}$ be such a set
that its image in $\ov{R}_{n,d}$ is a basis for
$\ov{R}_{n,d}(\De)$. For any multidegree
$\De=(\de_1,\ldots,\de_d)$ define $\MSG_{\De}$ by the following
way:
$$\MSG_{\De}=\MSG_{\de_{\si(1)}, \ldots,\de_{\si(d)}}|_{x_{ij}(r)\to
x_{ij}(\si(r)),\,i,j=\ov{1,n},\, r=\ov{1,d}},$$ %
where $\si\in S_d$, $\de_{\si(1)}\geq \cdots \geq\de_{\si(d)}$.
Then, the set $\MSG=\cup_{\de_1,\ldots,\de_d\geq 0} \MSG_{\De}$ is
a m.h.s.g. for $R_{n,d}$.
\end{remark}

Further, we assume that $n=3$, unless it is stated otherwise.

Let $B_{\De}$ be the basis for $\OK{\De}$ from
Proposition~\ref{prop_p0} and Theorems~\ref{theo2},~\ref{theo3}.
For $u\in K\LA F_d\RA^{\#}$ denote $\tr(u)= \tr(u|_{x_i\to X_i,\,
i=\ov{1,d}})\in R_{3,d}$.
\begin{theo}\label{theo4}
For  multidegree $\De=(\de_1,\ldots,\de_d)$, where $\de_1\geq
\cdots \geq\de_d$, $d\geq1$, define $\MSG_{\De}\subset R_{3,d}$:

1) the case of $p\neq3$:

\noindent if $d\geq2$ and $\de_d=1$, then
$\MSG_{\De}=\{\tr(ux_d)|\,u\in B_{(\de_1,\ldots,\de_{d-1})}\}$,

\noindent if $\De=2^3$, then
$\MSG_{\De}=\{\tr(X_1^2X_2^2X_3^2)\}$,

\noindent if $\De=2^2$, then $\MSG_{\De}=\{\tr(X_1^2X_2^2)\}$,

\noindent if $\De=3^2$, then
$\MSG_{\De}=\{\tr(X_1^2X_2^2X_1X_2)\}$,

\noindent if $d=1$, $\De=k$, $k=\ov{1,3}$, then
$\MSG_{\De}=\{\si_k(X_1)\}$,

\noindent for others $\De$ we define $\MSG_{\De}=\emptyset$;

2) the case of $p=3$:

\noindent if $d\geq2$ and $\de_d=1,2$, then
$\MSG_{\De}=\{\tr(ux_d^{\de_d})|\,u\in
B_{(\de_1,\ldots,\de_{d-1})}\}$,

\noindent if $\De=3^{2k}$, $k>0$, or $\De=3^{6k+1}$, $k>0$, then
$\MSG_{\De}=\{\tr(u)|\,u\in B_{\De}\}$,

\noindent if $d=1$, $\De=k$, $k=\ov{1,3}$, then
$\MSG_{\De}=\{\si_k(X_1)\}$,

\noindent for others $\De$ we define $\MSG_{\De}=\emptyset$.

Then, the set $\MSG$ from Remark~\ref{remark1} is a minimal system
of generators for $R_{3,d}$.
\end{theo}

In order to prove the Theorem we need some statements
from~\cite{Lopatin} and its corollaries:
\begin{lemma}\label{lemmaLop2}
1. Let $H\in A_{3,d-1}$. Then $\tr(HX_d)$ is decomposable if and
only if $\Phi(H)=0$ in $N_{3,d}$.

2. If $\tr(HX_d^2)\equiv0$, where $H\in A_{3,d-1}$, then
$\Phi(H)x_d+x_d\Phi(H)=0$ in $N_{3,d}$.

3. Let $p=3$, $H\in A_{3,d-1}$. Then $\tr(HX_d^2)\equiv0$ if and
only if $\Phi(H)=0$ in $N_{3,d-1}$.

4. If $\tr(u)$ is indecomposable, where $u$ is a word, then there
are canonical words $u_i$, $\mdeg(u)=\mdeg(u_i)$, and $\al_i\in
K$, such that $\tr(u)\equiv\sum\al_i \tr(u_i)$.

5. We have $\si_2(UV)\equiv\tr(U^2V^2)$, where $U,V\in A_{3,d}$.

6. Elements  $\si_2(X_1)$, $\det(X_1)$ are indecomposable.

7. If $p\neq3$, then $\tr(X_1^2X_2^2X_3^2) +
\tr(X_1^2X_3^2X_2^2)\equiv0$. For any $p$, the element
$\tr(X_1^2X_2^2X_3^2)$ is indecomposable.

8. Let $p=3$, $u_i\in F_d$ are words, $\mdeg(u_i)=\De$. Then
$\sum\nolimits_i \alpha_i\tr(u_i)\equiv0$ if and only if
$\sum_i\alpha_iu_i=0$ is a consequence of the system of identities
${\sys_{\De}}$ and identities $uv=vu$, where $u,v\in F_d^{\#}$ and
$\mdeg(uv)=\De$.

9. For $u,v\in F_d$, where $\mdeg(uv)=3^{2k}$, $k>0$, or
$\mdeg(uv)=3^{6k+1}$, $k>0$, we have $uv=vu$ in $N_{3,d}$.

10. Let $D$ be the explicit upper bound on degrees of elements of
a m.h.s.g. for $R_{3,d}$, where $d\geq2$. Then

if $p=0$ or $p>3$, then $D=6$;

if $p=2$, then $D=\left\{
\begin{array}{ccl}
d+2&,& d\geq4\\
6&,& d=2,3\\
\end{array}
\right.$;

if $p=3$ and $d=6k+r$, where $r\in\{3,5\}$, $k\geq0$, then
$D=3d-1$.

\end{lemma}
\begin{proof} All items except for $3$, $7$ are proven in~\cite{Lopatin}.

$3$. If $\tr(HX_d^2)\equiv0$, then  $\Phi(H)x_d+x_d\Phi(H)=0$ in
$N_{3,d}$ (see item~$2$). By item~$4$ of Lemma~\ref{lemma10} we
have $2\Phi(H)=0$ in $N_{3,d}$. The converse follows from
item~$1$.

$7$. Let $p\neq3$. The identity $x_1^2x_2x_3^2=0$ in $N_{3,d}$
(Lemma~\ref{lemmaLop} item~$4$) implies
$\tr(X_1^2X_2X_3^2X_2)\equiv0$ (see item~$1$). On the other hand,
the identity $x_2x_3^2x_2=-x_2^2x_3^2-x_3^2x_2^2$ in $N_{3,d}$
(see identity~\Ref{eq_t2}) implies $\tr(X_1^2X_2X_3^2X_2)\equiv
-\tr(X_1^2X_2^2X_3^2) - \tr(X_1^2X_3^2X_2^2)$ (see item~$1$). The
claim is proved.

Assuming $\tr(X_1^2X_2^2X_3^2)\equiv0$, we get
$x_1^2x_2^2x_3+x_3x_1^2x_2^2=0$ in $N_{3,d}$ by item~$2$. Thus
$x_1^2x_2^2x_1=-x_1x_1^2x_2^2=0$ in $N_{3,d}$; that is a
contradiction to item~$2$ of Lemma~\ref{lemmaLop}.
\end{proof}

\begin{proof_theo4}
By items~$4,5$ of Lemma~\ref{lemmaLop2}, we have that $R_{3,d}$ is
generated by $\{\si_2(X_i), \det(X_i), \tr(u)|\, u\in
F_d^{\#}\mbox{ is a canonical word},\,i=\ov{1,d}\}$.

Let $p\neq3$. The claim follows from items~$1$,~$6,$~$7$ and~$10$
of Lemma~\ref{lemmaLop2}.

Let $p=3$. The case of $\de_d=1,2$ follows from items~$1$,~$3$ of
Lemma~\ref{lemmaLop2}. The case of $\De=3^d$ follows from
items~$8$,~$9$,~$10$ of Lemma~\ref{lemmaLop2}.
\end{proof_theo4}

\begin{center} { ACKNOWLEDGEMENTS} \end{center}

The author is supported by RFFI (grant 04.01.00489). The author is
grateful to the referee whose comments considerably improved the
paper.



\end{document}